\documentclass{amsart} 
\usepackage{amssymb,amsmath,latexsym,times,xcolor,hyperref,tikz}%,enumitem}

\hypersetup{colorlinks=true, linkcolor=blue, citecolor=blue}

\numberwithin{equation}{section}

\theoremstyle{plain}
\newtheorem{thm}{Theorem}[section]

\newtheorem{cor}[thm]{Corollary}

\newtheorem{lemma}[thm]{Lemma}

\newtheorem*{lem*}{Lemma}
\newtheorem*{thm*}{Theorem}
\newtheorem*{cor*}{Corollary}

\theoremstyle{definition}

\newcommand{\del}{\backslash}

\DeclareMathOperator{\cl}{cl}
\DeclareMathOperator{\Tr}{Tr}
\DeclareMathOperator{\PG}{PG}
\DeclareMathOperator{\AG}{AG}
\DeclareMathOperator{\GF}{GF}

\title[Extremal matroids characterized by valuative
invariants]{Characterizations of certain matroids by maximizing
  valuative invariants} \author[J.~Bonin]{Joseph E.~Bonin} \address
{Department of Mathematics\\ The George Washington University\\
  Washington, D.C.\ 20052, USA} \email {jbonin@gwu.edu} \date{\today}


\begin{document}

\begin{abstract}
  Luis Ferroni and Alex Fink recently introduced a polytope of all
  unlabeled matroids of rank $r$ on $n$ elements, and they showed that
  the vertices of this polytope come from matroids that can be
  characterized by maximizing a sequence of valuative invariants.  We
  prove that a number of the matroids that they conjectured to yield
  vertices indeed do (these include cycle matroids of complete graphs,
  projective geometries, and Dowling geometries), and we give
  additional examples (including truncations of cycle matroids of
  complete graphs, Bose-Burton geometries, and binary and free spikes
  with tips).  We prove a special case of a conjecture of Ferroni and
  Fink by showing that direct sums of uniform matroids yield vertices
  of their polytope, and we prove a similar result for direct sums
  whose components are in certain restricted classes of matroids.
\end{abstract}

\maketitle

\section{Introduction}\label{sec:intro}

By a result of Derksen and Fink \cite{DerksenFink}, when the base
polytope of a matroid is written as an indicator function, the base
polytope of any matroid of rank $r$ on $n$ elements is a unique linear
combination of base polytopes of nested matroids of rank $r$ on $n$
elements and the coefficients are integers.  Matroid invariants assign
the same value to isomorphic matroids.  Valuative invariants are
matroid invariants $f$ that take values in an abelian group and for
which, for any matroid $M$, the image $f(M)$ can be found by taking
the integer linear combination of the base polytopes of nested
matroids that give the base polytope of $M$ and replacing each nested
matroid by its image under $f$.  (See Section \ref{ssec:VIandG} for
more precise statements.)  Valuative invariants include many
well-known invariants, such as the Tutte polynomial, but they
encompass far more than deletion-contraction invariants.  Ferroni and
Fink \cite{LuisAlex} introduced two polytopes of matroids; one is the
convex hull of the coordinate vectors of the base polytopes of
matroids of rank $r$ and size $n$ relative to the base polytopes of
nested matroids of rank $r$ and size $n$; the other, $\Omega_{r,n}$,
is a counterpart for unlabeled matroids (that is, isomorphism types of
matroids).  They showed that the vertices of $\Omega_{r,n}$ come from
matroids, or sets of matroids that give the same point in
$\Omega_{r,n}$, that can be characterized by sequences of maximization
(or, equivalently, minimization) problems of valuative invariants.  We
illustrate such characterizations, which are the focus of this paper,
with Theorem \ref{thm:Kn} below, which shows that, as Ferroni and Fink
conjectured, the cycle matroid $M(K_{r+1})$ of the complete graph
$K_{r+1}$ yields a vertex.
\begin{quote}
  \emph{Assume that $M$ is a simple matroid of rank $r$ on
    $\binom{r+1}{2}$ elements, that all lines have $2$ or $3$
    elements, that all planes have $3$, $4$, or $6$ elements, and
    that, for all $h$ with $4\leq h< r$, all rank-$h$ flats of $M$
    have at most $\binom{h+1}{2}$ elements.  Then the number of
    $3$-point lines in $M$ is at most $\binom{r+1}{3}$; also, if the
    number of $3$-point lines in $M$ is $\binom{r+1}{3}$, then $M$ is
    isomorphic to $M(K_{r+1})$.}
\end{quote}
To see how this conforms to the description of the vertices of
$\Omega_{r,\binom{r+1}{2}}$, we use the fact that the number of flats
of a given rank and size in a matroid is a valuative invariant, so
maximizing the number of flats of rank $1$ and size $1$ limits the
scope to simple matroids.  Then we can minimize the number of lines of
sizes other than two or three to say that all lines in the matroids of
interest have two or three points, and likewise we can then impose the
remaining conditions on the sizes of the flats of each rank.  Relative
to those conditions, we show that the maximum number of $3$-point
lines is $\binom{r+1}{3}$, and, for the final link, we show that the
maximum is achieved precisely by $M(K_{r+1})$.  Thus, a sequence of
optimization problems for valuative invariants hones in on
$M(K_{r+1})$, so $M(K_{r+1})$ gives rise to a vertex of
$\Omega_{r,\binom{r+1}{2}}$.

As that example suggests, these problems have the flavor of many that
have long been studied, but with a new twist, namely, we are limited
to optimizing valuative invariants.  Thus, above we can not assume at
the outset that the number of $3$-point lines is $\binom{r+1}{3}$; we
can do that only after showing that that value is the maximum among
matroids that satisfy the prior constraints.  Identifying vertices of
Ferroni and Fink's polytopes $\Omega_{r,n}$ thus opens up a new and
interesting chapter in extremal matroid theory.  Appropriately, they
call the matroids that yield vertices in their polytope extremal
matroids.

After reviewing the necessary background in Section
\ref{sec:background}, we show that several matroids that Ferroni and
Fink conjectured to be extremal indeed are: $M(K_{r+1})$ (Theorems
\ref{thm:Kn} and \ref{thm:Kn2}) and Dowling geometries $Q_r(G)$ based
on groups (or quasigroups, if $r=3$; Theorems \ref{thm:Dowling} and
\ref{thm:Dowling2}) are extremal.  They conjectured that projective
geometries are extremal; in Section \ref{sec:PMD} we prove that, more
broadly, all perfect matroid designs (i.e., matroids in which flats of
the same rank have the same size) are extremal; this applies to
projective geometries, affine geometries, their truncations, and much
more.  Among matroids that Ferroni and Fink did not mention, we show
that truncations of $M(K_{r+1})$ to rank at least five (Theorem
\ref{thm:truncateMKn}), Bose-Burton geometries (Theorems \ref{thm:BBg}
and \ref{thm:BBg2}), and binary and free spikes with tips (Theorem
\ref{thm:spikes}) are extremal.  Ferroni and Fink conjecture that
direct sums of extremal matroids are extremal; we do not have a proof
of that conjecture, but provide evidence for it by showing that direct
sums of uniform matroids are extremal (Theorem \ref{thm:dirsumunif})
and giving a general result (Theorem \ref{thm:dirsumfamilies}) that
shows that direct sums of cycle matroids of complete graphs are
extremal, and likewise for projective geometries of the same order,
affine geometries of the same order, Dowling geometries based on a
group that is the only group of its order, and certain spikes
(Corollary \ref{cor:extKnPGAGD}).

\section{Background}\label{sec:background}

For basic matroid theory, we recommend \cite{oxley}; we follow the
notation used there.  All matroids considered in this paper are
finite.  We use $\mathbb{N}$ for the set of positive integers,
$\mathbb{N}_0$ for the set of nonnegative integers, and $[n]$ for the
set $\{1,2,\ldots,n\}$.

\subsection{A few items about matroids}
A set $X$ in a matroid $M$ is \emph{cyclic} if it is a union of
circuits, or, equivalently, $M|X$ has no coloops.  We are most
interested in \emph{cyclic flats}, that is, flats that are cyclic.
The set of cyclic flats of $M$, which is denoted $\mathcal{Z}(M)$, is
a lattice under inclusion.  The least element of $\mathcal{Z}(M)$ is
$\cl_M(\emptyset)$, the set of loops of $M$, and the greatest element
is $E(M)-\cl_{M^*}(\emptyset)$ since $\cl_{M^*}(\emptyset)$ is the set
of coloops of $M$.  Routine arguments give the equality
$r(X)=\min\{r(F)+|X-F|\,:\,F\in\mathcal{Z}(M)\}$ for all
$X\subseteq E(M)$, so knowing $E(M)$, the cyclic flats of $M$, and the
rank of each cyclic flat determines $M$.

We use the term cover in the lattice of flats of a matroid $M$ and in
its lattice of cyclic flats $\mathcal{Z}(M)$.  In general, $y$ is a
\emph{cover} of $x$, or $y$ \emph{covers} $x$, in an ordered set if
$x<y$ and there is no $z$ with $x<z<y$.  If $Y$ covers $X$ in the
lattice of flats, then $r(Y)=r(X)+1$, but if $Y$ covers $X$ in
$\mathcal{Z}(M)$, then $r(Y)-r(X)$ can be any integer in
$[r(M)-r(X)]$.

The \emph{single-element extensions of $M$ by $e$} are the matroids
$N$ with $E(N)=E(M)\cup \{e\}$ and $N\del e=M$.  \emph{Modular cuts}
of $M$ are sets $\mathcal{M}$ of flats of $M$ that have two
properties: (i) if $F\in\mathcal{M}$ and $F'$ is a flat of $M$ with
$F\subseteq F'$, then $F'\in\mathcal{M}$, and (ii) if
$F,F'\in\mathcal{M}$ and $(F,F')$ is a \emph{modular pair}, that is,
$r(F)+r(F')=r(F\cup F')+r(F\cap F')$, then $F\cap F'\in\mathcal{M}$.
If $N$ is a single-element extension of $M$ by $e$, then the set of
flats $F$ of $M$ for which $r_N(F)=r_N(F\cup \{e\})$ is a modular cut.
Conversely, if $\mathcal{M}$ is a modular cut of $M$ and
$e\not\in E(M)$, then the following sets give the flats of a
single-element extension of $M$ by $e$: the sets $F\cup \{e\}$ with
$F\in\mathcal{M}$; the flats $F$ of $M$ for which
$F\not\in\mathcal{M}$; the sets $F\cup \{e\}$ for which $F$ is a flat
of $M$ and is neither in $\mathcal{M}$ nor covered by a flat in
$\mathcal{M}$.  The \emph{free extension of $M$ by $e$}, denoted
$M+e$, is the single-element extension that is given by the modular
cut $\mathcal{M}=\{E(M)\}$.  Thus, $E(M)\cup\{e\}$ is the unique
cyclic flat of $M+e$ that contains $e$.

Matroids are labeled: the ground set matters, and, for instance, which
sets are bases matters.  However, sometimes all we care about is the
isomorphism type of the matroid.  In that case, we call the
isomorphism class of a matroid, or a preferred representative of that
class, an \emph{unlabeled matroid}.  The next item illustrates this
with a relevant class of matroids.

A matroid $M$ is \emph{nested} if $\mathcal{Z}(M)$ is a chain, that
is, for any cyclic flats $F$ and $G$ of $M$, either $F\subseteq G$ or
$G\subseteq F$.  Nested matroids are known by many names, including
Schubert matroids.  Often all we care about is the isomorphism type of
a nested matroid.  Natural representatives of isomorphism classes of
nested matroids are matroids on $[n]$ in which the nonempty cyclic
flats are initial intervals $[s]$ of $[n]$.  The number of unlabeled
nested matroids of rank $r$ on $n$ elements is $\binom{n}{r}$, as we
now justify.  Consider sequences of matroids $M_0,M_1,\ldots,M_n$,
where $M_0$ is empty, and for each $i\in[n]$, we get $M_i$ by adding
the element $i$ to $M_{i-1}$ either (a) as a coloop or (b) via free
extension; then $M_n$ is a nested matroid on $[n]$, and if option (a)
is applied $r$ times, then $M_n$ has rank $r$.  Since $E(M)\cup\{e\}$
is the unique cyclic flat of $M+e$ that contains $e$, consistent with
what we said about a natural representative, the set of nonempty
cyclic flats of $M_n$ is
$$\{[s]\,:\, s, \text{ but not } s+1, \text{ is added by free
  extension}\}.$$ The number of elements in $[s]$ that are added as
coloops is the rank of $[s]$.  There are $\binom{n}{r}$ options for
which elements are added as coloops, so there are $\binom{n}{r}$
unlabeled rank-$r$ nested matroids on $n$ elements.  We can encode
$M_n$ as a lattice path $P=s_1s_2\ldots s_n$ where $s_i=N$ (a north
step) if $i$ is a coloop of $M_i$, otherwise $s_i=E$ (an east step);
if $r(M_n)=r$, then the lattice path $P$ goes from $(0,0)$ to
$(n-r,r)$.  For instance, the lattice path $NENEENNE$ encodes the
unlabeled rank-$4$ nested matroid on eight elements in which the
cyclic flats have the following size-rank pairs $(|X|,r(X))$: $(0,0)$,
$(2,1)$, $(5,2)$, and $(8,4)$; the corresponding sets in the natural
representative are $\emptyset$, $[2]$, $[5]$, and $[8]$.

Recall the Scum Theorem \cite[Theorem 3.3.1]{oxley}: if $N$ is a minor
of $M$, then there are subsets $X$ and $Y$ of $E(M)$ for which
$N=(M/Y)\del X$ and $r(N)=r(M/Y)$; also, if $N$ has no loops, then $Y$
can be taken to be a flat of $M$.  In particular, a matroid $M$ has no
minor isomorphic to $U_{2,q+2}$ if and only if, for each
rank-$(r(M)-2)$ flat $F$, at most $q+1$ hyperplanes of $M$ contain
$F$.

We will review the matroids that may be less familiar (Dowling
geometries, spikes, and Bose-Burton geometries) at the beginnings of
the sections in which they are treated.

\subsection{The Tutte polynomial}

The Tutte polynomial is the most well-known enumerative invariant in
matroid theory (see, e.g., \cite{Tutte,handbook}).  For a matroid $M$,
its \emph{Tutte polynomial} is defined to be
$$T(M;x,y) = \sum_{A\subseteq
  E(M)}(x-1)^{r(M)-r(A)}(y-1)^{|A|-r(A)},$$ which is a generating
function for the $2^{|E(M)|}$ pairs $(|A|,r(A))$ for
$A\subseteq E(M)$.  Among the most basic properties of the Tutte
polynomial are
\begin{itemize}
\item[(i)] the \emph{deletion/contraction formula}: if $e\in E(M)$ is
  neither a loop nor a coloop, then
  $T(M;x,y)=T(M\del e;x,y)+T(M/e;x,y)$, and
\item[(ii)] the \emph{direct sum factorization}:
  $T(M\oplus N;x,y)=T(M;x,y)T(N;x,y)$,
\end{itemize}
both of which are routine to prove.  With induction and these
properties, it is easy to prove that all coefficients in $T(M;x,y)$
are nonnegative integers and if $E(M)\ne\emptyset$, then $T(M;x,y)$
has no constant term.  The following well-known result of Brylawski
\cite{Tom} is also easy to prove by induction on $|E(M)|$ using
properties (i) and (ii), along with a key result by Tutte
\cite{TutteC}: for any connected matroid $M$ with $|E(M)|\geq 2$ and
any $e\in E(M)$, at least one of $M\del e$ and $M/e$ is connected.

\begin{lemma}
  A matroid $M$ is connected if and only if the coefficient of at
  least one of $x$ or $y$ in its Tutte polynomial $T(M;x,y)$ is
  nonzero.  
\end{lemma}

\begin{cor}
  A matroid $M$ has at least $k$ connected components if and only if
  all terms $x^iy^j$ in $T(M;x,y)$ with $i+j<k$ have zero coefficient.
\end{cor}

The Tutte polynomial is the universal deletion/contraction invariant,
as the following result by Oxley and Welsh \cite{univtutte} (extending
a result of Brylawski \cite{Tom}) makes precise.

\begin{thm}\label{thm:TutteUniversal}
  Let $R$ be a commutative ring with unity.  For any
  $u,v,\sigma,\tau\in R$ with $\sigma$ and $\tau$ nonzero, there is a
  unique function $t$ from the set of all matroids into $R$ that has
  the following properties.
  \begin{enumerate}
  \item[(1)] If $E(M)=\emptyset$, then $t(M) = 1$.
  \item[(2)] If $e$ is a loop of $M$, then
    $t(M) = v\cdot t(M\backslash e)$.
  \item[(3)] If $e$ is a coloop of $M$, then
    $t(M) = u\cdot t(M/e)$.
  \item[(4)] If $e$ is neither a loop nor a coloop, then
    $t(M) = \sigma\cdot t(M \backslash e) + \tau\cdot t(M/e)$.
  \end{enumerate}
  Furthermore, $t$ is the following evaluation of the Tutte polynomial
  $T(M;x,y)$:
  $$t(M)= \sigma^{|E(M)|-r(M)} \tau^{r(M)} T(M;u/\tau,v/\sigma).$$
\end{thm}

A matroid $M$ is \emph{Tutte unique} if the equality
$T(M;x,y)=T(N;x,y)$ holds only when $N$ is isomorphic to $M$.  Many
earlier characterizations of matroids by numerical invariants that we
will use, such as the characterization of $M(K_{r+1})$ in Theorem
\ref{thm:oldKnchar}, were proven to show that these matroids are Tutte
unique, but in most cases the characterizations do not fit the pattern
of maximizing a sequence of valuative invariants that is required in
order to show that the matroids yield vertices of the polytope of
unlabeled matroids.  Some Tutte-unique matroids do not yield vertices
of the polytope of unlabeled matroids (e.g., the unlabeled nested
matroid that is encoded by the lattice path $NENE$; see \cite[Examples
3.3 and 4.5]{LuisAlex}), and some matroids that yield vertices are not
Tutte unique (e.g., Dowling geometries of the same rank based on
nonisomorphic groups of the same order).  However, there are strong
connections between these two notions.  For more about Tutte-unique
matroids, see \cite{TUnTEq}.

\subsection{Matroid base polytopes, valuative invariants, and the
  $\mathcal{G}$-invariant}\label{ssec:VIandG}

We discuss matroid base polytopes very lightly; readers wanting more
on this topic are referred to \cite{AFR,DerksenFink,FS}.  To be able
to use simpler notation, when discussing matroid base polytopes it is
common to take the ground set of a matroid to be $[n]$.  Given a
matroid $M$ on $[n]$, a basis $B$ of $M$ can be encoded by its
\emph{characteristic vector} $\mathbf{v}_B$ in $\mathbb{R}^n$, that
is, the vector in which entry $i$ is $1$ if and only if $i\in B$,
otherwise entry $i$ is $0$.  The convex hull in $\mathbb{R}^n$ of the
vectors $\mathbf{v}_B$, as $B$ ranges over all bases of $M$, is the
\emph{matroid base polytope} $\mathcal{P}(M)$.

For a subset $X$ of $\mathbb{R}^n$, its \emph{indicator function}
$\bar{\mathbb{I}}(X):\mathbb{R}^n\to \mathbb{R}$ is given by
$$\bar{\mathbb{I}}(X)(x)=
\begin{cases}
  1, & \text{ if } x\in X,\\
  0, & \text{ if } x\not\in X.
\end{cases}$$
We focus on  indicator functions of base polytopes of  matroids.  We
call $\bar{\mathbb{I}}(\mathcal{P}(M))$ the  \emph{indicator
  function of} $M$.   The next  result is due to Derksen and
Fink \cite[Theorem 5.4]{DerksenFink}.

\begin{thm}\label{thm:dfnested}
  Let $M$ be a rank-$r$ matroid on $[n]$ and let $N_1,N_2,\ldots,N_t$
  be the nested matroids of rank $r$ on $[n]$.  There are unique
  integers $a_1,a_2,\ldots,a_t$ for which
  \begin{equation}\label{eq:lcnested}
    \bar{\mathbb{I}}(\mathcal{P}(M))=a_1\cdot
    \bar{\mathbb{I}}(\mathcal{P}(N_1))+ a_2\cdot
    \bar{\mathbb{I}}(\mathcal{P}(N_2))+\cdots+a_t\cdot
    \bar{\mathbb{I}}(\mathcal{P}(N_t)).
  \end{equation}
\end{thm}

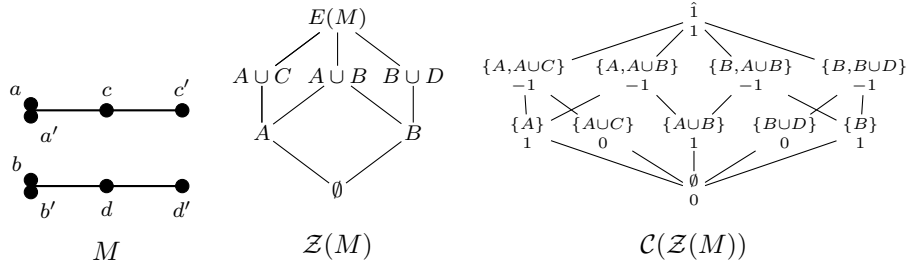
\begin{figure}
  \centering
  \begin{tikzpicture}[scale=1]
    \filldraw (0,4.92) node[below right] {\footnotesize $a'$} circle
    (2.5pt);%
    \filldraw (0,5.08) node[above left] {\footnotesize $a$} circle
    (2.5pt);%
    \filldraw (1,5) node[above=2pt] {\footnotesize $c$} circle
    (2.5pt);%
    \filldraw (2,5) node[above=2pt] {\footnotesize $c'$} circle
    (2.5pt);%
    \filldraw (0,3.92) node[below right] {\footnotesize $b'$} circle
    (2.5pt);%
    \filldraw (0,4.08) node[above left] {\footnotesize $b$} circle
    (2.5pt);%
    \filldraw (1,4) node[below=2pt] {\footnotesize $d$} circle
    (2.5pt);%
    \filldraw (2,4) node[below=2pt] {\footnotesize $d'$} circle
    (2.5pt);%
    \draw[thick](0,4)--(2,4);%
    \draw[thick](0,5)--(2,5);%
    \node at (1,3.15) {$M$};%
  \end{tikzpicture}
  \hspace{6pt}
  \begin{tikzpicture}[scale=1]
    \node[inner sep = 0.3mm] (em) at (1,0) {\footnotesize
      $\emptyset$};%
    \node[inner sep = 0.3mm] (s1) at (0,0.75) {\footnotesize $A$};%
    \node[inner sep = 0.3mm] (s3) at (2,0.75) {\footnotesize $B$};%
    \node[inner sep = 0.3mm] (s12) at (0,1.5) {\footnotesize
      $A\cup C$};%
    \node[inner sep = 0.3mm] (s13) at (1,1.5) {\footnotesize
      $A\cup B$};%
    \node[inner sep = 0.3mm] (s23) at (2,1.5) {\footnotesize
      $B\cup D$};%
    \node[inner sep = 0.3mm] (s123) at (1,2.25) {\footnotesize
      $E(M)$};%

    \foreach \from/\to in
    {em/s1,em/s3,s1/s12,s1/s13,s3/s23,s3/s13,s12/s123,
      s13/s123,s23/s123} \draw(\from)--(\to);%

    \foreach \from/\to in{s1/s12,s12/s123} \draw (\from)--(\to);%

    \node at (1,-0.75) {$\mathcal{Z}(M)$};%
  \end{tikzpicture}
  \hspace{6pt}
  \begin{tikzpicture}[scale=1]
    \node[inner sep = 0.3mm] (em) at (0,0) {\small
       $\genfrac{}{}{0pt}{}{\emptyset}{0}$};%

    \node[inner sep = 0.3mm] (s1) at (-2.2,0.75) {\small
      $\genfrac{}{}{0pt}{}{\{A\}}{1}$};%
    \node[inner sep = 0.3mm] (s2) at (-1.2,0.75) {\small
      $\genfrac{}{}{0pt}{}{\{A\cup C\}}{0}$};%
    \node[inner sep = 0.3mm] (s3) at (0,0.75) {\small
      $\genfrac{}{}{0pt}{}{\{A\cup B\}}{1}$};%
    \node[inner sep = 0.3mm] (s4) at (1.2,0.75) {\small
      $\genfrac{}{}{0pt}{}{\{B\cup D\}}{0}$};%
    \node[inner sep = 0.3mm] (s5) at (2.2,0.75) {\small
      $\genfrac{}{}{0pt}{}{\{B\}}{1}$};%

    \node[inner sep = 0.3mm] (t1) at (-2.25,1.5) {\small
      $\genfrac{}{}{0pt}{}{\{A,A\cup C\}}{-1}$} ;%
    \node[inner sep = 0.3mm] (t2) at (-0.75,1.5) {\small
      $\genfrac{}{}{0pt}{}{\{A,A\cup B\}}{-1}$};%
    \node[inner sep = 0.3mm] (t3) at (0.75,1.5) {\small
      $\genfrac{}{}{0pt}{}{\{B,A\cup B\}}{-1}$ };%
    \node[inner sep = 0.3mm] (t4) at (2.25,1.5) {\small
      $\genfrac{}{}{0pt}{}{\{B,B\cup D\}}{-1}$};%
    
    \node[inner sep = 0.3mm] (t) at (0,2.25) {\small
      $\genfrac{}{}{0pt}{}{\hat{1}}{1}$};%

    \foreach \from/\to in {em/s1,em/s2,em/s3,em/s4,em/s5,
      s1/t1,s1/t2,s2/t1,s3/t2,s3/t3,
      s4/t4,s5/t4,s5/t3,t/t1,t/t2,t/t3,t/t4} \draw(\from)--(\to);%

    \node at (0,-0.75) {$\mathcal{C}(\mathcal{Z}(M))$};%
\end{tikzpicture}  
\caption{A matroid $M$, its lattice $\mathcal{Z}(M)$ of cyclic flats,
  and the lattice $\mathcal{C}(\mathcal{Z}(M))$ of chains of cyclic
  flats that include $\emptyset$ and $E(M)$, with $\mu(C,\hat{1})$
  shown beneath $C$.  For readability, we suppress $\emptyset$ and
  $E(M)$ in the chains and let $A=\{a,a'\}$, $B=\{b,b'\}$,
  $C=\{c,c'\}$, and $D=\{d,d'\}$.}\label{fig:MuExample1}
\end{figure}

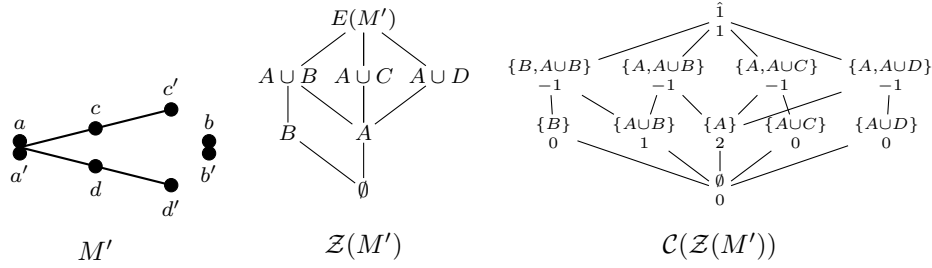
\begin{figure}
  \centering
  \begin{tikzpicture}[scale=1]
    \filldraw (5.75,4.42) node[below=1] {\footnotesize $a'$} circle
    (2.5pt);%
    \filldraw (5.75,4.58) node[above=1] {\footnotesize $a$} circle
    (2.5pt);%
    \filldraw (6.75,4.75) node[above=2pt] {\footnotesize $c$} circle
    (2.5pt);%
    \filldraw (7.75,5) node[above=2pt] {\footnotesize $c'$} circle
    (2.5pt);%
    \filldraw (6.75,4.25) node[below=2pt] {\footnotesize $d$} circle
    (2.5pt);%
    \filldraw (7.75,4) node[below=2pt] {\footnotesize $d'$} circle
    (2.5pt);%
    \filldraw (8.25,4.42) node[below =1] {\footnotesize $b'$} circle
    (2.5pt);%
    \filldraw (8.25,4.58) node[above =1] {\footnotesize $b$} circle
    (2.5pt);%
    \draw[thick](5.75,4.5)--(7.75,4);%
    \draw[thick](5.75,4.5)--(7.75,5);%
    \node at (6.75,3.15) {$M'$};%
  \end{tikzpicture}
  \hspace{6pt}
  \begin{tikzpicture}[scale=1]
    \node[inner sep = 0.3mm] (em) at (7,0) {\footnotesize
      $\emptyset$};%
    \node[inner sep = 0.3mm] (1) at (6,0.75) {\footnotesize
      $B$};%
    \node[inner sep = 0.3mm] (3) at (7,0.75) {\footnotesize
      $A$};%
    \node[inner sep = 0.3mm] (12) at (7,1.5) {\footnotesize
      $A\cup C$};%
    \node[inner sep = 0.3mm] (13) at (6,1.5) {\footnotesize
      $A\cup B$};%
    \node[inner sep = 0.3mm] (23) at (8,1.5) {\footnotesize
      $A\cup D$};%
    \node[inner sep = 0.3mm] (123) at (7,2.25) {\footnotesize
      $E(M')$};%

    \foreach \from/\to in
    {em/1,em/3,3/12,1/13,3/23,3/13,12/123,13/123,23/123}
    \draw(\from)--(\to);%

    \foreach \from/\to in {3/12,12/123} \draw (\from)--(\to);%

    \node at (7,-0.75) {$\mathcal{Z}(M')$};%
  \end{tikzpicture}
 \hspace{6pt}
  \begin{tikzpicture}[scale=1]
    \node[inner sep = 0.3mm] (em) at (0,0) {\small
       $\genfrac{}{}{0pt}{}{\emptyset}{0}$};%

    \node[inner sep = 0.3mm] (s1) at (-2.2,0.75) {\small
      $\genfrac{}{}{0pt}{}{\{B\}}{0}$};%
    \node[inner sep = 0.3mm] (s2) at (-1,0.75) {\small
      $\genfrac{}{}{0pt}{}{\{A\cup B\}}{1}$};%
    \node[inner sep = 0.3mm] (s3) at (0,0.75) {\small
      $\genfrac{}{}{0pt}{}{\{A\}}{2}$};%
    \node[inner sep = 0.3mm] (s4) at (1,0.75) {\small
      $\genfrac{}{}{0pt}{}{\{A\cup C\}}{0}$};%
    \node[inner sep = 0.3mm] (s5) at (2.2,0.75) {\small
      $\genfrac{}{}{0pt}{}{\{A\cup D\}}{0}$};%

    \node[inner sep = 0.3mm] (t1) at (-2.25,1.5) {\small
      $\genfrac{}{}{0pt}{}{\{B,A\cup B\}}{-1}$} ;%
    \node[inner sep = 0.3mm] (t2) at (-0.75,1.5) {\small
      $\genfrac{}{}{0pt}{}{\{A,A\cup B\}}{-1}$};%
    \node[inner sep = 0.3mm] (t3) at (0.75,1.5) {\small
      $\genfrac{}{}{0pt}{}{\{A,A\cup C\}}{-1}$ };%
    \node[inner sep = 0.3mm] (t4) at (2.25,1.5) {\small
      $\genfrac{}{}{0pt}{}{\{A,A\cup D\}}{-1}$};%
    
    \node[inner sep = 0.3mm] (t) at (0,2.25) {\small
      $\genfrac{}{}{0pt}{}{\hat{1}}{1}$};%

    \foreach \from/\to in {em/s1,em/s2,em/s3,em/s4,em/s5,
      s1/t1,s2/t2,s2/t1,s3/t2,s3/t3,s3/t4,
      s4/t3,s5/t4,t/t1,t/t2,t/t3,t/t4} \draw(\from)--(\to);%
    
    \node at (0,-0.75) {$\mathcal{C}(\mathcal{Z}(M'))$};%
  \end{tikzpicture}
  \caption{A second example of computing the coefficients
    $-\mu(C,\hat{1})$ in Equation (\ref{eq:lcnested}) following the
    same conventions as in Figure \ref{fig:MuExample1}.}\label{fig:MuExample2}
\end{figure}

An elegant and useful expression for the coefficients $a_i$ in
Equation (\ref{eq:lcnested}), due to Hampe \cite{Hampe}, uses the
M\"obius function.  (Readers who are not familiar with the M\"obius
function might consult Rota \cite{rota} or Stanley \cite{ec1}.  We use
the M\"obius function $\mu$ on lattices, and we need just two facts:
$\mu(x,x)=1$ for all $x$, and if $a<b$, then the sum of $\mu(x,b)$
over all $x$ with $a\leq x\leq b$ is $0$; with that, all values
$\mu(a,b)$ can be computed recursively.)  The set of cyclic flats of
any nested matroid $N_i$ for which $a_i\ne 0$ is a chain $C$ of cyclic
flats of $M$ that includes $\cl_M(\emptyset)$ and
$E(M)-\cl_{M^*}(\emptyset)$, the least and greatest cyclic flats of
$M$, and $r_{N_i}(X)=r_M(X)$ for all $X\in C$.  Let
$\mathcal{C}(\mathcal{Z}(M))$ consist of all chains in
$\mathcal{Z}(M)$ that include $\cl_M(\emptyset)$ and
$E(M)-\cl_{M^*}(\emptyset)$, ordered by inclusion, along with a
greatest element $\hat{1}$.  The intersection of two chains in
$\mathcal{C}(\mathcal{Z}(M))$ is in $\mathcal{C}(\mathcal{Z}(M))$, and
$\hat{1}$ is the greatest element, so $\mathcal{C}(\mathcal{Z}(M))$ is
a lattice.  The coefficient $a_i$ of
$\bar{\mathbb{I}}(\mathcal{P}(N_i))$, where $N_i$ is the nested
matroid whose cyclic flats and their ranks are those of the chain
$C\in \mathcal{C}(\mathcal{Z}(M))$, is the M\"obius value
$-\mu(C,\hat{1})$ in $\mathcal{C}(\mathcal{Z}(M))$.  Figures
\ref{fig:MuExample1} and \ref{fig:MuExample2} give examples.

The form of Equation (\ref{eq:lcnested}) motivates one of several
equivalent formulations of the notion of a valuation.  First, by a
\emph{matroid invariant}, we mean a function that is defined on the
set of matroids and that assigns the same value to isomorphic
matroids.  A matroid invariant $f$ taking values in an abelian group
$A$ is \emph{valuative}, or is a \emph{valuation}, if, for any
matroids $M_1,M_2,\ldots,M_t$ of rank $r$ on $[n]$ and any integers
$a_1,a_2,\ldots,a_t$, if 
$$a_1\cdot
\bar{\mathbb{I}}(\mathcal{P}(M_1))(x)+ a_2\cdot
\bar{\mathbb{I}}(\mathcal{P}(M_2))(x)+\cdots+a_t\cdot
\bar{\mathbb{I}}(\mathcal{P}(M_t))(x)=0$$ for all $x\in\mathbb{R}^n$,
then $a_1\cdot f(M_1)+a_2\cdot f(M_2)+\cdots+a_t\cdot f(M_t)$ is the
identity of $A$.  Thus, by Equation (\ref{eq:lcnested}), knowing the
values of a valuative invariant $f$ on (unlabeled) nested matroids
completely determines $f$.

For example, the indicator functions of the matroids in Figures
\ref{fig:MuExample1} and \ref{fig:MuExample2}, which must be
different, are written in terms of different nested matroids, and for
the matroid $M'$ in Figure \ref{fig:MuExample2}, we see a M\"obius
value, $-2$, that does not occur for the matroid $M$ in Figure
\ref{fig:MuExample1}.  However, any valuative invariant $f$ gives the
same values on these two matroids, namely,
$$4\cdot f(NENENEEE)-2\cdot f(NENNEEEE) -f(NNEENEEE),$$
where $f(P)$ is the value of $f$ at the unlabeled nested matroid that
is encoded by the lattice path $P$; the first term is from the four
coatoms in $\mathcal{C}(\mathcal{Z}(M))$ (or
$\mathcal{C}(\mathcal{Z}(M'))$) and the other terms are from the
atoms, some of which, besides $\emptyset$ and $E(M)$, include a cyclic
flat of rank $1$ and size $2$ (giving $NENNEEEE$) and the others
include a cyclic flat of rank $2$ and size $4$ (giving $NNEENEEE$).

Derksen and Fink \cite{DerksenFink} proved a counterpart of Theorem
\ref{thm:TutteUniversal} for valuative invariants; the universal
valuative invariant is the $\mathcal{G}$-invariant, to which we turn.
For a rank-$r$ matroid $M$ and permutation $\pi=e_1e_2\ldots e_n$ of
$E(M)$, the associated \emph{rank sequence} is the sequence
$\underline{r}(\pi) = r_1 r_2 \ldots r_n$ where $r_1 = r(\{e_1\})$ and
for $j \geq 2$,
$$r_j = r(\{e_1, e_2,\ldots,e_j\}) - r(\{e_1, e_2,\ldots,e_{j-1}\}).$$
Thus, $r$ entries of $\underline{r}(\pi)$ are $1$ and $n-r$ are $0$;
we call such sequences \emph{$(n,r)$-sequences}.  The rank sequence
$\underline{r}(\pi)$ shows the rank increases as one adjoins the
elements of $E(M)$, one at a time, to $\emptyset$ in the order given
by $\pi$.  The set $\{e_i\,:\,r_i=1\}$ is a basis of $M$.

Fix integers $n$ and $r$ with $0\leq r\leq n$.  For each
$(n,r)$-sequence $\underline{r}$, let $[\underline{r}]$ be a formal
symbol.  Let $\mathcal{G}(n,r)$ be the vector space, over a field $F$
of characteristic zero, that consists of all formal linear
combinations, with coefficients in $F$, of the symbols
$[\underline{r}]$.  For a matroid $M$ of rank $r$ on $n$ elements, its
\emph{$\mathcal{G}$-invariant} and its coefficients
$g_{\underline{r}}(M)$ are defined by
\begin{equation}\label{eq:GDef}
  \mathcal{G}(M) = \sum_{\underline{r}} g_{\underline{r}}(M)
  [\underline{r}] = \sum_{\pi} [\underline{r}(\pi)]
\end{equation}
where the second sum is over all $n!$ permutations $\pi$ of $E(M)$.
Thus, $g_{\underline{r}}(M)$ is the number of permutations $\pi$ of
$E(M)$ for which $\underline{r}(\pi)=\underline{r}$.  One can think of
$\mathcal{G}(M)$ as a generating function for rank sequences.

As mentioned above, Derksen and Fink \cite{DerksenFink} showed that
the $\mathcal{G}$-invariant is the universal valuative invariant: one
obtains any valuative invariant $f$ with values in an abelian group
$A$ by assigning values in $A$ to each symbol $[\underline{r}]$ and
extending linearly, that is, getting $f(M)$ by replacing each symbol
$[\underline{r}(\pi)]$ in Equation (\ref{eq:GDef}) by its image
$f([\underline{r}(\pi)])$ in $A$.  For example, for $A=\mathbb{Z}$
(the case of interest for this paper) and a fixed $(n,r)$-sequence
$\underline{r}$, by mapping $[\underline{r}]$ to $1$ and all other
rank sequences $[\underline{r}']$ to $0$, we see that the function
mapping $M$ to $g_{\underline{r}}(M)$ is valuative.  Also, integer
linear combinations of valuations with values in $\mathbb{Z}$ are
valuations.  Map each $[\underline{r}]$ to $[\underline{r}^*]$ where
we obtain $\underline{r}^*$ from $\underline{r}$ by replacing each $1$
by $0$, and each $0$ by $1$, and reversing the order of the entries;
it is routine to check that $\mathcal{G}(M)$ then maps to
$\mathcal{G}(M^*)$.  Thus, valuative invariants of $M^*$ are also
valuative invariants of $M$.

Speyer \cite{speyer} showed that the Tutte polynomial is a valuative
invariant.  From this and the definition of a valuation, it follows
easily that, for fixed nonnegative integers $i$ and $j$, the map that
sends a matroid $M$ to the coefficient of $x^iy^j$ in its Tutte
polynomial $T(M;x,y)$ is a valuative invariant.  The following theorem
collects these results along with valuative invariants that were
identified in \cite[Theorem 5.6]{catdata} and that we will use
heavily.

\begin{thm}\label{thm:listvalinv}
  Each of the following is a valuative invariant:
  \begin{enumerate}
  \item the map $M\mapsto \mathcal{G}(M)$,
  \item the map $M\mapsto \mathcal{G}(M^*)$,
  \item for any rank sequence $\underline{r}$, the map
    $M\mapsto g_{\underline{r}}(M)$,
  \item the map $M\mapsto T(M;x,y)$,
  \item for $i,j\in\mathbb{N}_0$, the map $M\mapsto t_{ij}$ where
    $t_{ij}$ is the coefficient of $x^iy^j$ in $T(M;x,y)$,
  \item the map $M\mapsto |\mathcal{I}(M)|$, where
    $\mathcal{I}(M)$ is the set of independent sets of $M$,
  \item for any $h,k,s_h,s_{h+1},\ldots,s_k\in \mathbb{N}_0$ for which
    $0\leq h\leq k\leq r(M)$, the map
    $M\mapsto F_{h,k}(M;s_h,\ldots,s_k)$, where
    $F_{h,k}(M;s_h,\ldots,s_k)$ is the number of flags
    $F_h\subsetneq F_{h+1}\subsetneq \cdots\subsetneq F_k$ of flats of
    $M$ with $r(F_i)=i$ and $ |F_i|=s_i$ for all $i$ with
    $h\leq i\leq k$,
  \item the map $M\mapsto f_k(M;s)$, where $f_k(M;s)$ is the number of
    rank-$k$ size-$s$ flats of $M$, and
  \item the map $M\mapsto f_k(M;s,c)$, where $f_k(M;s,c)$ is the
    number of rank-$k$ size-$s$ flats $X$ of $M$ for which $M|X$ has
    $c$ coloops.
  \end{enumerate}
\end{thm}

Item (8) is a special case of item (7).  By combining item (2) with
other items, it follows that the number of circuits of $M$ of a given
size is a valuative invariant, as is the number of cyclic sets of a
given size and nullity.  The number of cyclic flats of $M$ of a given
rank and size is a valuative invariant; set $c=0$ in item (9).  The
number of circuit-hyperplanes is a valuative invariant; in a rank-$r$
matroid $M$, it is $f_{r-1}(M;r,0)$.  Item (6) is a well-known
evaluation of the Tutte polynomial, as is the number of bases.

\subsection{Extremal matroids of the polytope of unlabeled matroids}

To associate analogs of matroid base polytopes to unlabeled matroids,
Ferroni and Fink use the \emph{symmetrized indicator function}
$\mathbb{I}(X):\mathbb{R}^n\to \mathbb{R}$ of a subset $X$ of
$\mathbb{R}^n$, which is given by
$$
\mathbb{I}(X)(x_1,x_2,\ldots,x_n)=\frac{1}{n!}\sum_{\sigma\in S_n}
\bar{\mathbb{I}}(X)(x_{\sigma(1)},x_{\sigma(2)},\ldots,x_{\sigma(n)})$$
where $S_n$ is the symmetric group on $[n]$. (The bar over
$\mathbb{I}$ is dropped in the symmetrized indicator function.) Note
that if $M$ and $N$ are isomorphic but unequal matroids on $[n]$, then
$\bar{\mathbb{I}}(\mathcal{P}(M))\ne \bar{\mathbb{I}}(\mathcal{P}(N))$
but $\mathbb{I}(\mathcal{P}(M))=\mathbb{I}(\mathcal{P}(N))$.  Thus,
letting the symmetrized indicator function of an unlabeled matroid be
the symmetrized indicator function of any (labeled) representative of
its isomorphism class makes sense.  Theorem \ref{thm:dfnested} carries
over to unlabeled matroids, but we must take into account that some
terms may merge.  For instance, for the matroid $M$ in Figure
\ref{fig:MuExample1}, for $\bar{\mathbb{I}}(\mathcal{P}(M))$, there
are two nested matroids whose cyclic flats are $\emptyset$, $E(M)$,
and a set of size two and rank $1$ (i.e., $A$ or $B$), and they give
different terms in Equation (\ref{eq:lcnested}), each with coefficient
$-1$; however, for $\mathbb{I}(\mathcal{P}(M))$, they give the same
term, namely, the symmetrized indicator function of the unlabeled
nested matroid that is encoded by the lattice path $NENNEEEE$, and its
coefficient is $-2$.  Theorem \ref{thm:dfnested} implies that the
symmetrized indicator function of any unlabeled rank-$r$ matroid on
$[n]$ can be written in exactly one way as an integer linear
combination of the symmetrized indicator functions of unlabeled
rank-$r$ nested matroids on $[n]$.  Fix a linear order on the set
$\mathcal{N}$ of all $\binom{n}{r}$ unlabeled rank-$r$ nested matroids
on $[n]$.  Let $U$ be the set of coordinate vectors of the symmetrized
indicator functions of unlabeled rank-$r$ matroids on $[n]$ relative
to those of the matroids in $\mathcal{N}$, so
$U\subsetneq\mathbb{R}^{\binom{n}{r}}$.  The convex hull of $U$ is the
\emph{polytope $\Omega_{r,n}$ of unlabeled rank-$r$ matroids on
  $[n]$}.  A matroid that gives a vertex in $\Omega_{r,n}$ is called
\emph{extremal}.  The starting point for our work is the following
result from \cite[Lemma 4.3]{LuisAlex}.

\begin{thm}\label{thm:vertbyinv}
  Fix integers $r$ and $n$ with $0\leq r\leq n$.  Let
  $f_1,f_2,\ldots,f_t$ be real-valued valuative invariants.  Let the
  set $S_0$ consist of one matroid of each isomorphism type among
  rank-$r$ matroids on $n$ elements, and, for each $i\in [t]$, set
  $$S_i=\{M\in S_{i-1}\,:\, f_i(M)\geq f_i(N) \text{ for all } N\in
  S_{i-1}\}.$$ If all matroids in $S_t$ yield the same point in the
  polytope $\Omega_{r,n}$, then that point is a vertex of
  $\Omega_{r,n}$, so all matroids in $S_t$ are then extremal.
\end{thm}

Theorem \ref{thm:vertbyinv} motivates characterizing certain matroids
by maximizing a sequence of real-valued valuative invariants.  Such
characterizations are the focus of this paper.  Note that $f$ is
valuative if and only if $-f$ is, so we can also minimize real-valued
valuative invariants.  In some settings maximizing is more natural,
while in others minimizing is more natural, so we refer to the generic
case as optimizing valuative invariants.

As noted at the beginning of \cite[Section 4.3]{LuisAlex}, a matroid
$M$ is extremal for $\Omega_{r,n}$ if and only if $M^*$ is extremal
for $\Omega_{n-r,n}$.  Valuative invariants of $M$ are valuative
invariants of $M^*$, so a characterization of $M$ by optimizing
valuative invariants also gives such a characterization of $M^*$.

We close this section with a useful observation that is implicit in
\cite{LuisAlex} and was made explicit by the referee, who also noted
its application to shortening some of the original proofs in this
paper.

\begin{lemma}\label{lem:refobs}
  For matroids $M$ and $N$ of rank $r$ on $n$ elements, if
  $\mathcal{G}(M)=\mathcal{G}(N)$, then $M$ and $N$ yield the same
  point in $\Omega_{r,n}$. Thus, $M$ is extremal if and only if $N$ is
  extremal.
\end{lemma}

Among the extremal matroids to which this lemma applies are perfect
matroid designs that have the same sequence of flat sizes (see the
next section; non-isomorphic projective planes of the same order are
familiar examples) and Dowling geometries of the same rank based on
non-isomorphic groups of the same order (see Section
\ref{sec:DL}). The lemma holds because if $M$ and $N$ yielded
different points $p_M$ and $p_N$ in $\Omega_{r,n}$, then there would
be a linear functional $\lambda$ on $\Omega_{r,n}$ for which
$\lambda(p_M)\ne\lambda(p_N)$, and that would give a valuative
invariant $f$ for which $f(M)\ne f(N)$, which contradicts having
$\mathcal{G}(M)=\mathcal{G}(N)$.

\section{Perfect matroid designs}\label{sec:PMD}

A \emph{perfect matroid design} is a matroid $M$ in which flats of the
same rank have the same size, that is, there are integers
$\alpha_0,\alpha_1,\ldots, \alpha_{r(M)}$ so that $|F|=\alpha_{r(F)}$
for all flats $F$ of $M$.  (See \cite{deza} and \cite{pmd}.)  Examples
include (finite) projective and affine geometries and their
truncations, as well as paving matroids of rank $r$ on $[n]$ for which
the hyperplanes are the blocks of a Steiner system $S(r-1,k,n)$ (i.e.,
a set of $k$-subsets (blocks) of $[n]$ for which each $(r-1)$-subset
of $[n]$ is a subset of exactly one block).

Nonisomorphic projective planes of order $q$ show that nonisomorphic
perfect matroid designs can have the same sequence of flat sizes.  In
contrast, any perfect matroid design that has the same sequence of
flat sizes as the projective geometry $\PG(r-1,q)$ with $r\geq 4$
satisfies the axioms for projective geometries and so is isomorphic to
$\PG(r-1,q)$, and likewise for the affine geometry $\AG(r-1,q)$.
(See, e.g., \cite{MK}.)

By \cite[Example 3.4]{catdata}, the $\mathcal{G}$-invariant of a
perfect matroid design $M$ is determined by its sequence
$\alpha_0,\alpha_1,\ldots, \alpha_{r(M)}$ of flat sizes.  Thus,
$f_i(M;\alpha_i)$, the number of rank-$i$ flats, is determined by
$\alpha_0,\alpha_1,\ldots, \alpha_{r(M)}$.  Also, since perfect
matroid designs with the same sequence of flat sizes have the same
$\mathcal{G}$-invariant, by Lemma~\ref{lem:refobs} they determine the
same point in the polytope of unlabeled matroids.  The valuative
invariant used below to show that perfect matroid designs are extremal
is defined in Theorem \ref{thm:listvalinv}.

\begin{thm}\label{thm:PGviaFlags}
  Perfect matroid designs are extremal.  If perfect matroid designs
  having the sequence of flat sizes
  $\alpha_0,\alpha_1,\ldots, \alpha_r$ exist, then they are precisely
  the rank-$r$ matroids on $\alpha_r$ elements that maximize the
  valuative invariant
  $N\mapsto F_{0,r}(N; \alpha_0,\alpha_1,\ldots, \alpha_r)$.  In
  particular, for $r\geq 4$ and a prime power $q$, among all simple
  rank-$r$ matroids on $\frac{q^r-1}{q-1}$ elements, $\PG(r-1,q)$ is
  the unique matroid that maximizes the valuative invariant
  $$N\mapsto
  F_{0,r}(N;0,1,q+1,q^2+q+1,\ldots,\frac{q^r-1}{q-1} ),$$ and among
  simple rank-$r$ matroids on $q^{r-1}$ elements, $\AG(r-1,q)$ is the
  unique matroid that maximizes the valuative invariant
  $N\mapsto F_{0,r}(N;0,1,q,q^2,\ldots,q^{r-1})$.
\end{thm}

\begin{proof}
  Recall that each element that is not in a flat $F$ is in exactly one
  flat that covers $F$.  Thus, if $0\leq i<r$ and $F$ is a rank-$i$
  flat of $M$ with $|F|=\alpha_i$, then at most
  $\frac{\alpha_r-\alpha_i}{\alpha_{i+1}-\alpha_i}$ flats cover $F$
  and have size $\alpha_{i+1}$, and there are
  $\frac{\alpha_r-\alpha_i}{\alpha_{i+1}-\alpha_i}$ such covers if and
  only if all covers of $F$ have size $\alpha_{i+1}$.  It follows that
  $$F_{0,r}(M;\alpha_0,\alpha_1,\ldots, \alpha_r) \leq
  \prod_{i=0}^{r-1}\frac{\alpha_r-\alpha_i}{\alpha_{i+1}-\alpha_i},$$
  with equality if and only if $M$ is a perfect matroid design with
  the sequence $\alpha_0,\alpha_1,\ldots, \alpha_r$ of flat sizes.  It
  now follows that perfect matroid designs are extremal since, as
  noted above, all perfect matroid designs with the sequence
  $\alpha_0,\alpha_1,\ldots, \alpha_r$ of flat sizes determine the
  same point in the polytope of matroids.
\end{proof}

\section{Bose-Burton geometries}

Bose-Burton geometries are matroids of the form
$\PG(r-1,q)\del \PG(k-1,q)$, that is, the matroid obtained by deleting
from $\PG(r-1,q)$ one of its rank-$k$ flats.  The next theorem, from
Bose and Burton \cite{bose}, is a counterpart, for projective
geometries, of Tur\'an's theorem; just as Tur\'an's theorem treats
subgraphs of $K_n$ that do not have another particular complete graph
as a subgraph, this theorem treats restrictions of $\PG(r-1,q)$ that
do not have another particular projective geometry over $\GF(q)$ as a
restriction.

\begin{thm}\label{thm:bose}
  If $M$ is a restriction of $\PG(r-1,q)$ that has no restriction
  isomorphic to $\PG(m-1,q)$, then
  $$|E(M)|\leq \frac{q^r - q^{r-m+1}}{q-1}.$$  
  Furthermore, equality holds if and only if $M$ is isomorphic to
  $\PG(r-1,q) \del \PG(r-m,q)$.
\end{thm}

That is one of many attractive properties of Bose-Burton geometries.
If $k=r-1$, then $\PG(r-1,q)\del \PG(k-1,q)$ is the affine geometry
$\AG(r-1,q)$, which is extremal since it is a perfect matroid design.
By the result below from \cite{Rachelle}, all other Bose-Burton
geometries except perhaps $\PG(r-1,q)\del \PG(r-3,q)$ are extremal
since conditions (1)--(3) can be cast as optimizing valuative
invariants.

\begin{thm}\label{thm:char}
  Let $q$ be a prime power.  Fix $r,k\in\mathbb{N}$ with $r\geq 4$ and
  $1\leq k \leq r-3$.  Let $M$ be a simple rank-$r$ matroid on
  $(q^r - q^k)/(q-1)$ elements that satisfies conditions
  \emph{(1)--(3)}.
  \begin{itemize}
  \item[(1)] If $F$ is a rank-$(r-k+1)$ flat of $M$, then
    $|F|<(q^{r-k+1} - 1)/(q-1)$.
  \item[(2)] For any hyperplane $H$ of $M$, either
    $$|H|=\frac{q^{r-1} - q^k}{q-1} \qquad \text{ or }\qquad
    |H|=\frac{q^{r-1} - q^{k-1}}{q-1}.$$ 
  \item[(3)] If $F$ is a rank-$(r-2)$ flat of $M$, then $|F|$ is one
    of
    $$\frac{q^{r-2} - q^k}{q-1}, \qquad \frac{q^{r-2} - q^{k-1}}{q-1}, \quad
    \text{ or, if } k>1,\quad \frac{q^{r-2} - q^{k-2}}{q-1}.$$
  \end{itemize}
  Then $M$ is isomorphic to $\PG(r-1,q)\del \PG(k-1,q)$.
\end{thm}

Theorems \ref{thm:BBg} and \ref{thm:BBg2} treat the remaining case: we
show that $\PG(r-1,q)\del \PG(r-3,q)$ is extremal.  We address $r=3$
and $r\geq 4$ separately, in part because there can be more than one
projective plane of order $q$, so we must show that all single-element
deletions of all projective planes of order $q$ give the same point in
the polytope of unlabeled matroids.

\begin{thm}\label{thm:BBg}
  Fix an integer $q\geq 2$.  Let $M$ be a rank-$3$ simple matroid with
  $q^2+q$ elements in which each line has $q$ or $q+1$ points.  The
  number of $q$-point lines in $M$ is at least $q+1$, and it is $q+1$
  if and only if $M=N\del x$ for some projective plane $N$ of order
  $q$ and element $x\in E(N)$.  All single-element deletions of
  projective planes of order $q$ are extremal.
\end{thm}

\begin{proof}
  A single-element deletion $N\del x$ of a projective plane $N$ of
  order $q$ has $q+1$ lines that each contain $q$ points, namely,
  $L-\{x\}$ as $L$ runs through the lines of $N$ with $x\in L$.

  Now assume that each line of a rank-$3$ simple matroid $M$ on
  $q^2+q$ elements has size $q$ or $q+1$.  If all lines that contain
  $e\in E(M)$ had $q+1$ points, then $|E(M)-\{e\}|=q^2+q-1$ would be
  divisible by $q$, but it is not.  Thus, each element is in at least
  one $q$-point line.  Along with the factorization $q^2+q=q(q+1)$,
  that implies that at least $q+1$ lines have $q$ points each, as
  claimed.  That factorization also implies that $M$ has $q+1$ lines
  with $q$ points each if and only if each element of $M$ is in one
  such line.  In that case the $q$-point lines must be pairwise
  disjoint and so, along with $E(M)$, they form a modular cut.  The
  corresponding single-element extension, say $N$, has $q^2+q+1$
  points and all lines have $q+1$ points, from which it follows that
  $N$ is a projective plane of order $q$.

  There may be nonisomorphic projective planes of order $q$ if
  $q\geq 9$, so to complete the proof that their single-element
  deletions are extremal, we show that, for $q\geq 9$, all
  single-element deletions of all projective planes of order $q$ give
  the same point in $\Omega_{3,q^2+q}$.  For a single-element deletion
  $M$ of a projective plane of order $q$, its indicator function
  $\bar{\mathbb{I}}(\mathcal{P}(M))$ depends on $M$, but, as we show
  next, from $q$ alone we get the relevant unlabeled nested matroids,
  their multiplicities, and their coefficients, so the symmetrized
  indicator function $\mathbb{I}(\mathcal{P}(M))$, and hence the
  corresponding point in $\Omega_{3,q^2+q}$, depends only on $q$, not
  on $M$.  The lattice of cyclic flats of $M$ consists of $\emptyset$,
  $q+1$ lines that each have $q$ points, $q^2$ lines that each have
  $q+1$ points, and $E(M)$.  Thus, the lattice
  $\mathcal{C}(\mathcal{Z}(M))$ has, besides $\hat{1}$,
  \begin{itemize}
  \item $q+1$ elements $C=\{\emptyset,L,E(M)\}$ covered by $\hat{1}$
    where $r(L)=2$ and $|L|=q$; for each, $\mu(C,\hat{1})=-1$,
  \item $q^2$ elements $C=\{\emptyset,L,E(M)\}$ covered by $\hat{1}$
    where $r(L)=2$ and $|L|=q+1$; for each, $\mu(C,\hat{1})=-1$,
  \item the least chain, $C=\{\emptyset,E(M)\}$; given the values
    above, including $\mu(\hat{1},\hat{1})=1$, we have
    $\mu(C,\hat{1})=q^2+q$.
  \end{itemize}
  Thus, all single-element deletions of projective planes of order $q$
  give the same point in $\Omega_{3,q^2+q}$.
\end{proof}

The argument in the last paragraph of the proof above also shows that
all single-element deletions of projective planes of order $q$ have
the same $\mathcal{G}$-invariant.  We are not aware of a prior
appearance of that result in the literature, so we could not appeal to
Lemma~\ref{lem:refobs}.

We treat $r\geq 4$ in Theorem \ref{thm:BBg2} below.  The proof uses
the result below from \cite{Uq}.

\begin{thm}\label{thm:rep}
  Let $r\geq 4$ be an integer and let $q$ be a prime power.  Any
  simple rank-$r$ matroid that has no minor isomorphic to $U_{2,q+2}$,
  the $(q+2)$-point line, and that has at least $q^{r-1}$ points is
  representable over $\GF(q)$.
\end{thm}

\begin{thm}\label{thm:BBg2}
  Fix a prime power $q$.  Assume that $M$ is a simple matroid of rank
  $r\geq 4$ on $q^{r-1}+q^{r-2}$ elements that satisfies the following
  conditions:
  \begin{itemize}  
  \item[(1)] if $2\leq i<r$ and $A$ is a rank-$i$ flat of $M$, then
    $|A|\in\{q^{i-1},q^{i-1}+q^{i-2}\}$, and
  \item[(2)] $F_{r-2,r-1}(M;q^{r-3}+q^{r-4}, q^{r-2})=0$, that is, no
    rank-$(r-2)$ flat with $q^{r-3}+q^{r-4}$ elements is a subset of a
    hyperplane with $q^{r-2}$ elements.
  \end{itemize}
  Then the minimum number of $q^{r-2}$-point hyperplanes in $M$ is
  $q+1$.  Furthermore, $M$ has $q+1$ such hyperplanes if and only if
  each element of $M$ is in one such hyperplane, which is equivalent
  to $M$ being isomorphic to $\PG(r-1,q)\del \PG(r-3,q)$.  Thus, for
  any $r\geq 4$, the matroid $\PG(r-1,q)\del \PG(r-3,q)$ is extremal.
\end{thm}

\begin{proof} We first show that $\PG(r-1,q)\del X$, where $X$ is a
  rank-$(r-2)$ flat of $\PG(r-1,q)$, has the properties above.  For
  any flat $F$ of $\PG(r-1,q)$ where $F\not\subseteq X$, if
  $r(F\cup X)=r-1$, then $r(F\cap X)=r(F)-1$, and otherwise
  $r(F\cap X)=r(F)-2$; these two cases yield the options for $|F-X|$
  in condition (1).  Now let $F$ be a rank-$(r-2)$ flat of
  $\PG(r-1,q)$ for which $|F-X|=q^{r-3} +q^{r-4}$, so
  $r(X\cap F)=r-4$, and so $r(X\cup F)=r$.  Thus, if $H$ is a
  hyperplane of $\PG(r-1,q)$ and $F\subsetneq H$, then
  $X\not\subseteq H$, so $r(H\cap X)=r-3$, and so
  $|H-X|=q^{r-2}+q^{r-3}$; thus, condition (2) holds.  The
  $q^{r-2}$-point hyperplanes of $\PG(r-1,q)\del X$ are the sets $H-X$
  where $H$ ranges over the hyperplanes of $\PG(r-1,q)$ that contain
  $X$, so there are $q+1$ such hyperplanes.

  Now assume that $M$ satisfies conditions (1) and (2).  To prove the
  lower bound of $q+1$ on the number of $q^{r-2}$-point hyperplanes,
  note that if, for some $i\geq 1$, the element $e$ is in a rank-$i$
  flat $F$ with $q^{i-1}$ elements and $i<r-1$, then having all
  rank-$(i+1)$ flats that contain $F$ have $q^i+q^{i-1}$ elements
  would yield an integer solution $k$ of the equation
  $q^{r-1}+q^{r-2}= q^{i-1}+kq^i$, and so of
  $q^{r-i}+q^{r-i-1}= 1+kq$, which has no integer solution; thus, $e$
  is in some rank-$(i+1)$ flat with $q^i$ elements.  Thus, by
  induction, each $e\in E(M)$ is in at least one $q^{r-2}$-point
  hyperplane. The lower bound on the number of $q^{r-2}$-point
  hyperplanes now follows since the fewest number of sets of size
  $q^{r-2}$ whose union is $E(M)$ is $q+1$ since
  $q^{r-1}+q^{r-2}= q^{r-2}(q+1)$.  Since each element is in at least
  one $q^{r-2}$-point hyperplane, that factorization also shows that
  $M$ has $q+1$ such hyperplanes if and only if each element of $M$ is
  in one such hyperplane.
  
  Now assume that $M$ has $q+1$ hyperplanes of size $q^{r-2}$, so each
  element of $M$ is in one such hyperplane.  We next use Theorem
  \ref{thm:rep} to show that $M$ is representable over $\GF(q)$.  To
  show that $M$ has no $U_{2,q+2}$-minor, by the Scum Theorem it
  suffices to show that each rank-$(r-2)$-flat $F$ of $M$ is contained
  in at most $q+1$ hyperplanes.  Each element is in one
  $q^{r-2}$-point hyperplane, so at most one such hyperplane contains
  $F$.  First assume that $|F|= q^{r-3}$ and that $i$ hyperplanes that
  contain $F$ have $q^{r-2}$ elements, so $i\in\{0,1\}$.  The number
  of $(q^{r-2}+q^{r-3})$-point hyperplanes that contain $F$ is then
  $$\frac{q^{r-1}+q^{r-2} -q^{r-3}-i(q^{r-2} -q^{r-3})}{q^{r-2}}.$$
  Canceling common factors shows that $q$ divides $q^2+q -1-i(q -1)$,
  so $q$ divides $i-1$.  Thus, $i=1$ and the displayed quantity
  simplifies to $q$, so $F$ is in $q+1$ hyperplanes, as needed.  Now
  assume that $|F|= q^{r-3}+q^{r-4}$.  By conditions (1) and (2), each
  hyperplane that contains $F$ has size $q^{r-2}+q^{r-3}$, so the
  number of such hyperplanes is
  $$\frac{q^{r-1}+q^{r-2} -
    (q^{r-3}+q^{r-4})}{q^{r-2}+q^{r-3}-(q^{r-3}+q^{r-4})}= q+1.$$
  Thus, by the Scum Theorem and Theorem \ref{thm:rep}, $M$ is
  representable over $\GF(q)$ and so is isomorphic to a restriction of
  $\PG(r-1,q)$.
  
  Finally, condition (1) implies that $M$ has no flat $F$ for which
  $M|F$ is isomorphic to $\PG(2,q)$, so, by Theorem \ref{thm:bose},
  $M$ is isomorphic to $\PG(r-1,q)\del \PG(r-3,q)$.  Since, for fixed
  $q$ and $r$, optimizing the sequence of valuative invariants that
  correspond to conditions above yields one matroid,
  $\PG(r-1,q)\del \PG(r-3,q)$, it is extremal.
\end{proof}

\section{Cycle matroids of complete graphs, and their truncations}\label{sec:Kr+1}

To prove that $M(K_{r+1})$ is extremal, we use the characterization of
$M(K_{r+1})$ below from \cite[Theorem 3.2]{bill}, which does not have
the form set out in Theorem \ref{thm:vertbyinv} since the hypotheses
include the numbers of $3$-point lines and $6$-point planes without
establishing that they are the greatest numbers of such flats.

\begin{thm}\label{thm:oldKnchar}
  The cycle matroid $M(K_{r+1})$ is the unique rank-$r$ simple matroid
  on ${r+1\choose 2}$ elements that has
  \begin{itemize}
  \item [(1)] ${r+1\choose 3}$ lines with three points,
  \item [(2)] no $5$-point planes, ${r+1\choose 4}$ planes with six
    points, no planes with more than six points, and
  \item [(3)] no rank-$4$ flats with more than ten points.
  \end{itemize}
\end{thm}

As discussed in Section \ref{sec:intro}, all criteria in the next
result optimize valuative invariants, so this shows that $M(K_{r+1})$
is extremal, as conjectured in \cite{LuisAlex}.  This result sets the
pattern for later results in the paper in that since it is immediate
that the matroid of interest, $M(K_{r+1})$ in this case, satisfies the
conditions in the result, we do not explicitly state or prove that.

\begin{thm}\label{thm:Kn}
  Assume that $M$ is a simple rank-$r$ matroid on $\binom{r+1}{2}$
  elements, that all lines have $2$ or $3$ elements, that all planes
  have $3$, $4$, or $6$ elements, and that, for all $h$ with
  $4\leq h<r$, all rank-$h$ flats of $M$ have at most $\binom{h+1}{2}$
  elements.  Then at most $\binom{r+1}{3}$ lines of $M$ have three
  elements, and if $\binom{r+1}{3}$ lines of $M$ have three elements,
  then $M$ is $M(K_{r+1})$.  Thus, $M(K_{r+1})$ is extremal.
\end{thm}

\begin{proof}
  In this proof, we call $3$-point lines \emph{large lines}.  The
  statement clearly holds for $r\leq 2$, so now assume that $r\geq 3$.

  Let $F$ be a rank-$h$ flat of $M$ where $3\leq h\leq r$.  Let an
  element $e\in F$ be in $m$ large lines in $F$, say $L_i$, where
  $i\in[m]$.  We claim that $m\leq h-1$.  If $h=3$, then
  $1+2m\leq |F|\leq 6$, so $m\leq 2$, as needed.  Now consider $h=4$.
  No three large lines are coplanar by the case $h=3$.  Since no plane
  of $M$ has five elements, each plane $\cl(L_i\cup L_j)$ contains one
  point $a_{i,j}$ that is in neither $L_i$ nor $L_j$.  The $1+2m$
  elements in the lines $L_j$, for $j\in [m]$, and the $m-1$ distinct
  points $a_{i,m}$, for $i\in[m-1]$, give $3m \leq |F| \leq 10$, so
  $m\leq 3$, as needed.  Now consider $h>4$.  As above, each
  difference $\cl(L_i\cup L_j)-(L_i\cup L_j)$ is a singleton, say
  $\{a_{i,j}\}$.  We claim that, apart from the obvious equality
  $a_{i,j}=a_{j,i}$, all $a_{i,j}$ are distinct.  If
  $a_{i,j}=a_{i,k}$, then the plane $\cl(L_i\cup \{a_{i,j}\})$ would
  contain $L_i$, $L_j$, and $L_k$, contrary to the case $h=3$.  If
  $a_{i,j}=a_{s,t}$ with $\{i,j\}\cap \{s,t\}=\emptyset$, then the
  rank-$4$ flat $\cl(L_i\cup L_s\cup \{a_{i,j}\})$ would contain
  $L_i$, $L_j$, $L_s$, and $L_t$, contrary to the case $h=4$.  Since
  all $a_{i,j}$ are distinct, the rank-$h$ flat $F$ contains at least
  $1+2m+\binom{m}{2}$ elements, so
  $1+2m+\binom{m}{2}\leq \binom{h+1}{2}$, that is,
  $\binom{m+2}{2}\leq \binom{h+1}{2}$, so $m\leq h-1$, as claimed.
  
  Let $M$ have $k$ large lines.  The number of pairs $(L,e)$, where
  $L$ is a large line and $e\in L$, is $3k$.  There are
  $\binom{r+1}{2}$ points $e$ and each is in at most $r-1$ large
  lines, so $3k\leq \binom{r+1}{2}(r-1)$, that is,
  $k\leq\binom{r+1}{3}$, as the bound in the theorem asserts.

  Now assume that $M$ has $\binom{r+1}{3}$ large lines.  The counting
  arguments above show that each element $e$ of $M$ is in $r-1$ large
  lines, any $h-1$ of which span a flat of rank $h$.  Also, $E(M)$ is
  the union of the large lines $L_1,L_2,\ldots,L_{r-1}$ that contain
  $e$ along with the $\binom{r-1}{2}$ singletons
  $\{a_{i,j}\}=\cl(L_i\cup L_j)-(L_i\cup L_j)$ for distinct
  $i,j\in[r-1]$.  This notation applies throughout the next three
  paragraphs, which deduce more about the $6$-point planes of $M$.

  We claim that, for each $i\in[r-1]$, the planes $\cl(L_i\cup L_j)$,
  for $j\in[r-1]-\{i\}$, account for all $6$-point planes that contain
  $L_i$.  Since $\cl(L_i\cup L_j)= \cl(L_i\cup \{a_{i,j}\})$, each
  other plane that contains $L_i$ is $\cl(L_i\cup \{a_{s,t}\})$ for
  some $\{s,t\}\subseteq[r-1]-\{i\}$.  Since
  $\cl(L_i\cup \{a_{s,t}\})$ is in the rank-$4$ flat
  $\cl(L_i\cup L_s\cup L_t)$, and the two $6$-point planes
  $\cl(L_i\cup L_s)$ and $\cl(L_i\cup L_t)$ contain nine of the ten
  points in $\cl(L_i\cup L_s\cup L_t)$, the plane
  $\cl(L_i\cup \{a_{s,t}\})$ is $L_i\cup \{a_{s,t}\}$ and so has only
  four points, which proves the claim.

  Similarly, rank-$4$ flats that contain $L_i\cup L_j$ besides those
  of the form $\cl(L_i\cup L_j\cup L_k)$ have the form
  $F=\cl(L_i\cup L_j\cup\{a_{s,t}\})$ for some
  $\{s,t\}\subseteq [r-1]-\{i,j\}$.  Now $F$ is contained in the
  rank-$5$ flat $F'=\cl(L_i\cup L_j\cup L_s\cup L_t)$, and
  $\cl(L_i\cup L_j\cup L_s)$ and $\cl(L_i\cup L_j\cup L_t)$ contain
  $14$ of the $15$ points in $F'$, so $|F|=7$.  Thus, for all
  $\{i,j\}\subseteq[r-1]$, the flats $\cl(L_i\cup L_j\cup L_k),$ for
  $k\in[r-1]-\{i,j\}$, are the only $10$-element rank-$4$ flats that
  contain $L_i\cup L_j$.

  We claim that if $P$ is a $6$-point plane with
  $\{e,a_{i,j}\} \subsetneq P$, then $P=\cl(L_i\cup L_j)$.  If not,
  then $P\cap\cl(L_i\cup L_j)=\{e,a_{i,j}\}$.  Set
  $F=P\cup \cl(L_i\cup L_j)$.  Now $|F|=10$ and $r(F)=4$, so $F$ is a
  flat, and so $F=\cl(L_i\cup L_j\cup L_k)$, for some
  $k\in[r-1]-\{i,j\}$, by the previous paragraph.  Thus,
  $L_k\subsetneq P$ and $L_h\not\subseteq P$ for all
  $h\in[r-1]-\{k\}$, contrary to what we proved two paragraphs above.
  Thus, the $6$-point planes that contain $e$ are the planes
  $\cl(L_i\cup L_j)$, for $\{i,j\}\subseteq[r-2]$.  Any point in a
  $6$-point plane can be $e$, so each such plane contains four large
  lines.

  Let $k$ be the number of $6$-point planes of $M$.  Consider the
  triples $(P,L,e)$ where $P$ is a $6$-point plane, $L$ is a large
  line in $P$, and $e\in L$.  The number of such triples is
  $k\cdot 4\cdot 3=12k$ by the conclusion above.  By focusing on $e$
  first, and then $L$, and finally $P$, we see that the number of
  triples is $\binom{r+1}{2}(r-1)(r-2)$, where $r-2$ accounts for
  choosing a second line containing $e$ to determine the plane $P$.
  Thus, $12k=\binom{r+1}{2}(r-1)(r-2)$, so $k=\binom{r+1}{4}$.

  This shows that the criteria in Theorem \ref{thm:oldKnchar} hold, so
  $M$ is $M(K_{r+1})$.
\end{proof}

The next result gives two more characterizations of $M(K_{r+1})$ by
optimizing valuative invariants.  Each of statements (a) and (b) below
shows that $M(K_{r+1})$ is extremal.  The second suits our work in
Section \ref{sec:dirsum} on direct sums.  The proof uses Theorem
\ref{thm:Kn}.

\begin{thm}\label{thm:Kn2}
  Let $M$ be a simple matroid of rank $r\geq 2$ on $\binom{r+1}{2}$
  elements.  Assume that for all $i\in[r-1]$ and rank-$i$ flats $F$,
  either $|F|=\binom{i+1}{2}$ or $|F|\leq \binom{i}{2}+1$.  Call a
  rank-$i$ flat $F$ \emph{large} if $|F|=\binom{i+1}{2}$, and call a
  flag
  $F_0\subsetneq F_1\subsetneq F_2\subsetneq \cdots \subsetneq F_r$ of
  flats \emph{large} if each $F_i$ is large. 
  \begin{itemize}
  \item[(a)] The number of large hyperplanes is at most $r+1$, and if
    it is $r+1$, then $M$ is $M(K_{r+1})$.
  \item[(b)] The number of large flags is at most $(r+1)!/2$, and if
    it is $(r+1)!/2$, then $M$ is $M(K_{r+1})$.
  \end{itemize}
\end{thm}

\begin{proof}
  We prove the bound in (a) by induction on $r$.  An inequality that
  we use in the induction step holds for $r\geq 5$, so we check all
  cases $r<5$ first.  If $r=2$, then $M$ is $U_{2,3}$, so the
  assertion holds.  For $r=3$, since $|E(M)|=6$, each point is in at
  most two large lines; thus, if there are $k$ large lines, then the
  number of pairs $(e,L)$ where $L$ is a large line and $e\in L$ is,
  on the one hand, $3k$, and, on the other, at most $6\cdot 2$, so
  $3k\leq 12$, that is, $k\leq 4$, as needed.

  Now consider $r=4$, so $|E(M)|=10$.  Note that the intersection of
  any two large planes $P$ and $P'$ is a line with two or three
  elements.  First assume that $P\cap P'=\{a,b\}$.  Thus,
  $P\cup P'=E(M)$, so any other large plane $Q$ must intersect each
  of $P$ and $P'$ in a large line; let $Q\cap P=L$ and
  $Q\cap P'=L'$.  Now $Q=L\cup L'$, so $L\cap L'=\emptyset$, and
  so $\{a,b\}\cap(L\cup L')=\emptyset$.  However, only one large line
  in $P$ can be disjoint from $\{a,b\}$, and similarly for $P'$, so
  there is at most one large plane besides $P$ and $P'$, consistent
  with the bound.  Now assume that any two large planes of $M$
  intersect in a large line.  By counting, a large line is in at most
  two large planes.  Since any large plane $P$ has at most four large
  lines and each other large plane must intersect $P$ in one of those
  lines, and each such line is in at most one large plane besides $P$,
  there are at most five large planes, as needed.  It is relevant for
  the last step of the proof to note that if $M$ has five large
  planes, then any two of them intersect in a large line.
  
  Now assume that $r\geq 5$ and that the bound on the number of large
  hyperplanes holds in the case of rank $r-1$.  Let
  $H_1,H_2,\ldots,H_t$ be the large hyperplanes of $M$.  For any
  $\{i,j\}\subseteq [t]$, we have
  $|H_i\cup H_j|\leq |E(M)|=\binom{r+1}{2}$, so
  \begin{align*}
    |H_i\cap H_j|
    &\,= |H_i|+|H_j|-|H_i\cup H_j|\\
    &\,\geq 2 \binom{r}{2}- \binom{r+1}{2}\\
    &\,=  \binom{r-1}{2}- 1\\
    &\, > \binom{r-2}{2}+ 1
  \end{align*}
  where the last line holds since $r\geq 5$.  Therefore $H_i\cap H_j$
  is a large rank-$(r-2)$ flat.  Since $H_i-(H_i\cap H_j)$ and
  $H_j-(H_i\cap H_j)$ each contain $\binom{r}{2}-\binom{r-1}{2}=r-1$
  out of the $\binom{r+1}{2}-\binom{r-1}{2}=2r-1$ elements in
  $E(M)-(H_i\cap H_j)$, it follows that $H_i$ and $H_j$ are the only
  large hyperplanes that contain $H_i\cap H_j$, so all $\binom{t}{2}$
  intersections $H_i\cap H_j$, with $\{i,j\}\subseteq[t]$, are
  different.  Thus, $H_i\cap H_t$, for $i\in[t-1]$, are $t-1$
  different large hyperplanes of $M|H_t$, so the induction assumption
  applied to $M|H_t$ gives $t-1 \leq r$, so $t\leq r+1$, which
  completes the proof of the bound stated in (a).

  The bound of $(r+1)!/2$ on the number of large flags
  $F_0\subsetneq F_1\subsetneq \cdots \subsetneq F_r$ now follows
  since there are at most $r+1$ options for $F_{r-1}$, and then at
  most $r$ options for $F_{r-2}$, and so on, ending with at most three
  options for $F_1$.

  The argument shows that in order for $M$ to have $(r+1)!/2$ large
  flags, it must have $r+1$ large hyperplanes, so proving the second
  part of (a) will also complete the proof of (b).  Let
  $H_1,H_2,\ldots,H_{r+1}$ be the large hyperplanes of $M$.  As we
  showed, each $H_i$ contains $r$ large rank-$(r-2)$ flats, namely,
  $H_i\cap H_j$ for $j\in[r+1]-\{i\}$.  With that, a routine induction
  shows that if $I\subsetneq[r+1]$, then $\bigcap_{i\in I}H_i$ is a
  large rank-$(r-|I|)$ flat.  Thus, the only large hyperplanes that
  contain $\bigcap_{i\in I}H_i$ are $H_i$ with $i\in I$, so all such
  intersections are different.  Therefore the intersections of $r-2$
  large hyperplanes account for $\binom{r+1}{r-2}=\binom{r+1}{3}$
  different large lines of $M$, and so $M$ is $M(K_{r+1})$ by Theorem
  \ref{thm:Kn}.
\end{proof}

We now turn to truncations of $M(K_{r+1})$ of rank at least $5$.  In
\cite{AnnaMarc}, de Mier and Noy showed that these truncations are
Tutte unique.  For that, they first developed a numerical
characterization of these matroids similar to that of $M(K_{r+1})$ in
Theorem \ref{thm:oldKnchar}, adding an assumption on the number of
independent sets.  We characterize truncations of $M(K_{r+1})$ of rank
at least $5$ as optimizing certain valuative invariants and so show
that these matroids are extremal.  The proof below does not rest on
prior characterizations.  Theorem \ref{thm:truncateMKn} gives yet
another characterization of $M(K_{r+1})$ (set $k=0$), although it is
limited to $r\geq 5$.

Let $\Tr^k(N)$, where $0\leq k\leq r(N)$, denote the matroid obtained
from $N$ by truncating $k$ times.  Thus, $\Tr^k(N)$ has rank $r(N)-k$,
and its independent sets are those of $N$ that have size at most
$r(N)-k$.

\begin{thm}\label{thm:truncateMKn}
  Fix $r,k\in\mathbb{N}\cup \{0\}$ with $r-k\geq 5$.  Let $M$ be a
  simple rank-$(r-k)$ matroid on $\binom{r+1}{2}$ elements in which
  $|F|\leq \binom{r(F)+1}{2}$ for all flats $F\ne E(M)$ and no planes
  have five elements.  Then at most $\binom{r+1}{3}$ lines of $M$ have
  three elements.  If $\binom{r+1}{3}$ lines of $M$ have three
  elements, then $M$ has at most as many independent sets as
  $\Tr^k(M(K_{r+1}))$.  If, in addition, $M$ has the same number of
  independent sets as $\Tr^k(M(K_{r+1}))$, then $M$ is
  $\Tr^k(M(K_{r+1}))$.  Thus, $\Tr^k(M(K_{r+1}))$ is extremal.
\end{thm}

\begin{proof}
  As earlier, we call $3$-point lines \emph{large lines}.  The
  argument in the second paragraph of the proof of Theorem
  \ref{thm:Kn} applies and shows that for any $e\in E(M)$ and
  $h\in [r-k-1]$, the element $e$ is in at most $h-1$ large lines in
  any rank-$h$ flat.  Let $e$ be in $i$ large lines of $M$, say
  $L_1,L_2,\ldots,L_i$.  Now $|L_1\cup L_2\cup \cdots\cup L_i|=2i+1$,
  and, by the argument in the proof of Theorem \ref{thm:Kn}, the sets
  $\cl(L_s\cup L_t)-(L_s\cup L_t)$, for $1\leq s<t\leq i$, are
  disjoint singletons.  Thus,
  $1+2i+\binom{i}{2}\leq |E(M)|=\binom{r+1}{2}$, from which we get
  $i\leq r-1$.  With that, the argument in the third paragraph in the
  proof of Theorem \ref{thm:Kn} applies and shows that $M$ has at most
  $\binom{r+1}{3}$ large lines, as asserted in the theorem.

  Now assume that $M$ has $\binom{r+1}{3}$ large lines.  As in the
  proof of Theorem \ref{thm:Kn}, it follows that each element $e$ of
  $M$ is in $r-1$ large lines, any $h-1$ of which span a rank-$h$ flat
  if $h< r-k$.  Furthermore, $E(M)$ is the union of the $r-1$ large
  lines $L_1,L_2,\ldots,L_{r-1}$ that contain $e$ along with the
  $\binom{r-1}{2}$ singletons
  $\{a_{i,j}\}=\cl(L_i\cup L_j)-(L_i\cup L_j)$ for
  $\{i,j\}\subsetneq[r-1]$.  As in the proof of Theorem \ref{thm:Kn},
  for each $i\in[r-1]$, the planes $\cl(L_i\cup L_j)$, for
  $j\in[r-1]-\{i\}$, account for all $6$-point planes that contain
  $L_i$.  Consider $\{i,j\}\subsetneq[r-1]$ and $f\in L_i-\{e\}$.
  Since $f$ is in a large line in the plane $\cl(L_i\cup L_j)$, our
  previous conclusion, applied to $f$ and the large lines that contain
  it, shows that $f$ is in two large lines in the plane
  $\cl(L_i\cup L_j)$.  From this, it follows that each $6$-point plane
  of $M$ that contains at least one large line contains four large
  lines.

  We next identify the elements and large lines of $M$ with those of
  $M(K_{r+1})$.  Fix an element of $M$, label it $a_{r,r+1}$, and let
  $L_1,L_2,\ldots,L_{r-1}$ be the large lines that contain
  $a_{r,r+1}$.  For $i$ and $j$ with $1\leq i<j< r$, let the singleton
  $\cl(L_i\cup L_j)-(L_i\cup L_j)$ be $\{a_{i,j}\}$.  Arbitrarily
  label the elements in $L_1-\{ a_{r,r+1}\}$ as $a_{1,r}$ and
  $a_{1,r+1}$.  For each $i$ with $2\leq i <r$, let
  $L_i\cap\cl(\{a_{1,r},a_{1,i}\})=\{a_{i,r}\}$ and
  $L_i\cap\cl(\{a_{1,r+1},a_{1,i}\})=\{a_{i,r+1}\}$, so
  $L_i=\{a_{r,r+1},a_{i,r},a_{i,r+1}\}$.  To show that this labeling
  agrees with the large lines of $M(K_{r+1})$, where $K_{r+1}$ has
  vertices $v_1,v_2,\ldots,v_{r+1}$ and we want $a_{i,j}$ to
  correspond to the edge $\{v_i,v_j\}$, we want to show that
  \begin{itemize}
  \item[(a)] for $1<i<j<r$, both $\{a_{i,r},a_{j,r},a_{i,j}\}$ and
    $\{a_{i,r+1},a_{j,r+1},a_{i,j}\}$ are lines, 
  \item[(b)] for $1\leq i<j<k<r$, the set
    $\{a_{i,j},a_{i,k},a_{j,k}\}$ is a line,
  \end{itemize}
  and that we have identified all large lines of $M$.  For the set
  $\{a_{i,r},a_{j,r},a_{i,j}\}$ in assertion (a), note that the plane
  $\cl(\{a_{1,r},a_{i,r},a_{j,r}\})$ is contained in the rank-$4$ flat
  $\cl(L_1\cup L_i\cup L_j)$ and it contains $a_{1,i}$ and $a_{1,j}$
  but none of $L_1$, $L_i$, and $L_j$.  Therefore $a_{i,j}$ must be
  the sixth element of the plane $\cl(\{a_{1,r},a_{i,r},a_{j,r}\})$.
  Thus, $\{a_{i,r},a_{j,r},a_{i,j}\}$ is the intersection of the
  planes $\cl(\{a_{1,r},a_{i,r},a_{j,r}\})$ and $\cl(L_i\cup J_j)$,
  and so is a line.  A similar argument shows that
  $\{a_{i,r+1},a_{j,r+1},a_{i,j}\}$ is a line.  For (b), by (a) and
  the labeling, the plane $P_r=\cl(\{a_{i,j},a_{j,k},a_{j,r}\})$
  contains both $a_{i,r}$ (which is in $\cl(\{a_{i,j},a_{j,r}\})$) and
  $a_{k,r}$ (which is in $\cl(\{a_{j,k},a_{j,r}\})$), and so $P_r$
  contains $a_{i,k}$, which is in $\cl(\{a_{i,r},a_{k,r}\})$;
  likewise, the plane $P_{r+1}=\cl(\{a_{i,j},a_{j,k},a_{j,r+1}\})$
  contains both $a_{i,r+1}$ and $a_{k,r+1}$, and so it contains
  $a_{i,k}$.  Thus, $\{a_{i,j},a_{i,k},a_{j,k}\}$ is $P_r\cap P_{r+1}$
  and so is a line.  Finally, the $r-1$ large lines
  $L_1,L_2,\ldots,L_{r-1}$, the $2\binom{r-1}{2}$ other large lines in
  the planes $\cl(L_i\cup L_j)$, and the $\binom{r-1}{3}$ large lines
  identified in (b) account for all $\binom{r+1}{3}$ lines of $M$.
  With this identification, we can now assume that $M$ and
  $\Tr^k(M(K_{r+1}))$ have the same ground set and the same $3$-point
  lines.

  We next show that sets that are independent in $M$ are independent
  in $\Tr^k(M(K_{r+1}))$; equivalently, every circuit $C$ of
  $\Tr^k(M(K_{r+1}))$ is dependent in $M$.  This holds if $C$ spans
  $\Tr^k(M(K_{r+1}))$ since then $|C|>r(M)$.  It also holds if $|C|=3$
  by our work above.  Now assume that $4\leq|C|\leq r-k$.  In
  $K_{r+1}$, the set $C$ contains the edges of a cycle; fix an edge
  $e$ in $C$, let $v$ run through the vertices in that cycle that are
  not incident with $e$, and consider the $3$-cycles determined by $e$
  and $v$, which are large lines of both $M$ and $\Tr^k(M(K_{r+1}))$
  that contain $e$. Each edge of $C$ is in either one of these lines
  or the closure of two of them, so, in both $M$ and
  $\Tr^k(M(K_{r+1}))$, the set $C$ is in the closure of the union of
  $|C|-2$ large lines that contain $e$.  That union has rank $|C|-1$
  in $M$, so $C$ is dependent in $M$.

  By what we just showed, $M$ has at most as many independent sets as
  $\Tr^k(M(K_{r+1}))$; furthermore, if $M$ and $\Tr^k(M(K_{r+1}))$
  have the same number of independent sets, then they have the same
  independent sets and so are the same matroid.
\end{proof}

\section{Dowling geometries}\label{sec:DL}

We first recall the idea behind Dowling geometries.  Let
$\{b_1,b_2,\ldots,b_r\}$ be a basis of the projective geometry
$\PG(r-1,q)$ over the finite field $\GF(q)$.  Let $Q_r(\GF^*(q))$ be
the restriction of $\PG(r-1,q)$ to the set
$\bigcup_{\{i,j\}\subseteq[r]}\cl(\{b_i,b_j\})$.  The independent
sets, circuits, flats, etc., in $Q_r(\GF^*(q))$ can be written using
only the group $\GF^*(q)$ of units of $\GF(q)$.  Dowling
\cite{Dowling} showed how to replace $\GF^*(q)$ by an arbitrary finite
group.

There are many ways to define Dowling geometries, including via
group-labeled partial partitions (see \cite{Dowling}) and
group-labeled graphs (see \cite{oxley}).  Since cyclic flats play a
big role in the polytope of matroids, we define Dowling geometries by
giving their cyclic flats and the rank of each cyclic flat.  In a
matroid $M$ with no loops, for each connected component $X$ of a
nonempty cyclic flat, $|X|>1$ and $X$ is a \emph{connected flat},
i.e., $X$ is a flat and $M|X$ is connected.  We give the connected
flats $X$ of Dowling geometries with $|X|>1$, and use them to give the
cyclic flats.

\begin{figure}
  \centering
  \begin{tikzpicture}[scale=0.85]
    \draw[thick] (-30:2) -- (90:2) -- (210:2) -- (-30:2); %

    \filldraw (120:1.15) node[left=1] {\footnotesize$-1_{(1,2)}$}
    circle (3pt);%
    \filldraw (90:2) node[left=1] {\footnotesize$b_1$} circle (3pt);%
    \filldraw (180:1.15) node[left=1] {\footnotesize$1_{(1,2)}$}
    circle (3pt);%
    
    \filldraw (240:1.15) node[below=1] {\footnotesize$-1_{(2,3)} \,\,\,\, $}
    circle (3pt);%
    \filldraw (210:2) node[left=1] {\footnotesize$b_2$} circle (3pt);%
    \filldraw (300:1.15) node[below=2] {\footnotesize$\,\,\,\, 1_{(2,3)}$}
    circle (3pt);%

    \filldraw (60:1.15) node[right=1] {\footnotesize$-1_{(1,3)}$}
    circle (3pt);%
    \filldraw (330:2) node[right=1] {\footnotesize$b_3$} circle
    (3pt);%
    \filldraw (0:1.15) node[right=1] {\footnotesize$1_{(1,3)}$} circle
    (3pt);%

    \draw[thick] (180:1.15) to [out=270,in=180] (240:1.15)  to
    [out=0,in=270] (60:1.15) ;%

    \draw[thick] (180:1.15) to [out=300,in=180] (300:1.15)  to
    [out=0,in=270] (0:1.15) ;%

    \draw[very thick,dashed] (120:1.15) to [out=250,in=180] (240:1.15) to
    [out=0,in=240] (0:1.15) ;%

    \draw[very thin] (120:1.15) to [out=250,in=180] (240:1.15) to
    [out=0,in=240] (0:1.15) ;%

    \draw[very thick,dotted] (120:1.15) to [out=280,in=180] (300:1.15)
    to [out=0,in=280] (60:1.15) ;%

    \draw[very thin] (120:1.15) to [out=280,in=180] (300:1.15) to
    [out=0,in=280] (60:1.15) ;%
    
  \end{tikzpicture}
  \caption{The Dowling geometry $Q_3(G)$ where $G$ is the sign group,
    $\{1,-1\}$.} 
  \label{fig:Q3sign}
\end{figure}
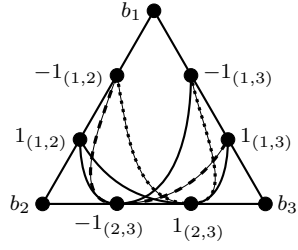

Let $G$ be a finite group, written with multiplicative notation. Fix
$r\geq 3$.  The ground set of the \emph{rank-$r$ Dowling geometry
  $Q_r(G)$ over $G$} is
$$\{b_1, b_2,\ldots, b_r\}\cup\{ a_{(i,j)}\,:\, a \in G \text{ and } 1 \leq
i < j \leq r\}.$$ Hence $Q_r(G)$ has $r+|G|\binom{r}{2}$ elements,
among which are $\binom{r}{2}$ copies of each element of $G$, each one
indexed by an ordered pair $(i,j)$ with $1\leq i<j\leq r$.  For $i>j$,
let $a_{(i,j)} = (a^{-1})_{(j,i)}$, where $a^{-1}$ is the inverse of
$a$ in $G$.  The connected flats $F$ with $|F|>1$ are of two types:
\begin{itemize}
\item[(L)] for $I\subseteq [r]$ with $|I|\geq 2$, the \emph{large
    flat}
  $$\{b_i\,:\,i\in I\}\cup \{ a_{(i,j)}\,:\, a \in G \text{ and }
  \{i,j\} \subseteq I\}$$ has rank $|I|$ and \emph{support} $I$;
\item[(S)] for a subset $J=\{j_1<j_2<\cdots<j_t\}$ of $[r]$ with
  $t\geq 3$ and (not necessarily distinct) $a_2,a_3,\ldots,a_t\in G$,
  the \emph{small flat}
  $$\{(a_h)_{(j_1,j_h)}\,:\,2\leq h\leq t\}\cup  \{
  (a_h^{-1}a_i)_{(j_h,j_i)}\,:\, 2\leq h<i\leq t\}$$ has rank $t-1$
  and \emph{support} $J$.
\end{itemize}
Thus, a large flat with support $I$ has size $|I|+|G|\binom{|I|}{2}$,
and a small flat with support $J$ has size $\binom{|J|}{2}$.  The
cyclic flats of $Q_r(G)$ are the unions of at most one large flat and
any number of small flats for which all supports are pairwise
disjoint; the rank of a cyclic flat is the sum of the ranks of the
large and small flats that it contains.  Figure \ref{fig:Q3sign}
illustrates the construction with $Q_3(G)$ where $G$ is the sign
group, $\{1,-1\}$.  In that example, the set $I=\{1,2\}$ gives the
large line $\{b_1,b_2,1_{(1,2)},-1_{(1,2)}\}$, and the set $J=[3]$ and
the group elements $a_2=a_3=-1$ give the dotted line
$\{-1_{(1,2)},1_{(2,3)},-1_{(1,3)}\}$.

It is easy to check that small flats are the sets $F$ for which, for
some $J\subseteq[r]$ with $|J|\geq 3$, (i) $F$ consists of one element
from each of the $\binom{|J|}{2}$ sets $\{a_{(i,j)}\,:\,a\in G\}$ for
$\{i,j\}\subsetneq J$, and (ii) if $a_{(i,j)},a'_{(j,k)}\in F$, where
$\{i,j,k\}\subseteq J$, then $(aa')_{(i,k)}\in F$.  Also, the
restriction of $Q_r(G)$ to a small rank-$h$ flat is isomorphic to
$M(K_{h+1})$ and the restriction to a large rank-$h$ flat is
isomorphic to $Q_h(G)$.  If $|G|=1$, then $Q_r(G)$ is isomorphic to
$M(K_{r+1})$, which we treated in Section \ref{sec:Kr+1}.  The results
in this section apply to Dowling geometries $Q_r(G)$ with $|G|>1$.

We can define $Q_3(G)$ when $G$ is a quasigroup (a Latin square), that
is, $G$ has a binary operation, and the equations $ax=b$ and $ya=b$
have unique solutions for all $a,b \in G$.  Quasigroups may lack an
identity and inverses, so the identification
$a_{(i,j)} = (a^{-1})_{(j,i)}$ mentioned above for groups makes no
sense; thus, we write $a_{(i,j)}$ only when $i<j$.  The small lines of
$Q_3(G)$ are $\{a_{(1,2)},a'_{(2,3)},(aa')_{(1,3)}\}$ for $a,a'\in G$.
For $r=3$, our results apply to $Q_3(G)$ where $G$ is any quasigroup
with $|G|>1$.

We will use the following characterization of Dowling geometries.  The
case of $r\geq 4$ was proven in \cite{DLChar1}; the extension to
$r=3$, along with a different proof, appeared in \cite{DLChar2}, which
also gave a geometric criterion for distinguishing rank-$3$ Dowling
geometries based on groups from those based on other quasigroups.
That criterion cannot be cast in terms of valuative invariants since
all Dowling geometries of the same rank based on groups (or, if $r=3$,
quasigroups) of the same order have the same $\mathcal{G}$-invariant
\cite[Proposition 4.8]{catdata}.

\begin{thm}\label{thm:axioms}
  A simple matroid $M$ of rank $r \geq 3$ is a Dowling geometry (based
  on a quasigroup if $r=3$) if and only if $M$ has elements
  $b_1, b_2,\dots,b_r$ for which
  \begin{enumerate}
  \item $E(M)$ is the union of the $\binom{r}{2}$ lines
    $\cl(\{b_i,b_j\})$, for $\{i,j\}\subsetneq[r]$,
  \item  $|\cl(\{b_i,b_j\})|>2$ for all $\{i,j\}\subsetneq[r]$, and
  \item if $\{i,j,k\}\subseteq [r]$,
    $e\in \cl(\{b_i,b_j\}) - \{b_i,b_j\}$, and
    $f\in \cl(\{b_i,b_k\}) - \{b_i,b_k\}$, then $|\cl(\{e, f\})|>2$.
  \end{enumerate}
\end{thm}

The set $\{b_1, b_2,\ldots, b_r\}$ is a basis of $M$.  When the group
or quasigroup is nontrivial, only one basis $\{b_1, b_2,\ldots, b_r\}$
satisfies properties (1)--(3) above.

The next result characterizes Dowling geometries of rank $r\geq 3$
other than $M(K_{r+1})$ by optimizing a sequence of valuative
invariants.  Since all Dowling geometries of rank $r$ based on groups
(or quasigroups if $r=3$) of order $g$ have the same
$\mathcal{G}$-invariant, the characterization and
Lemma~\ref{lem:refobs} show that they are extremal.

\begin{thm}\label{thm:Dowling}
  Fix integers $r\geq 3$ and $g\geq 2$.  Let $M$ be a simple matroid
  of rank $r$ with $|E(M)|=r+g\binom{r}{2}$ in which, for all
  $h\in[r-1]$, each rank-$h$ flat has at most $h+g\binom{h}{2}$
  elements.  Call a rank-$h$ flat $F$ \emph{large} if
  $|F|= h+g\binom{h}{2}$.  Assume that for any rank-$h$ flat $F$ that
  is not large, if $2<h<r$, then $|F|\leq h+g\binom{h-1}{2}$, and if
  $h=2$, then $|F|\leq 3$.  Then the following conclusions hold.
  \begin{itemize}
  \item[(a)] The matroid $M$ has at most $r$ large hyperplanes.
  \item[(b)] If $M$ has $r$ large hyperplanes, then $M$ has a basis
    $B$ for which $\cl(X)$ is large for all $X\subseteq B$.
  \item[(c)] If $M$ has $r$ large hyperplanes, then the greatest
    number of $3$-point lines in $M$ is $\binom{r}{3}g^2$, and any
    such matroid $M$ with $\binom{r}{3}g^2$ three-point lines is a
    rank-$r$ Dowling geometry based on a group (or quasigroup if
    $r=3$) of order $g$.
  \end{itemize}
  Thus, rank-$r$ Dowling geometries based on groups of order $g$ are
  extremal.
\end{thm}

\begin{proof}
  We first prove (a) by induction on $r$.  In addition to the base
  case $r=3$, we treat $r=4$ separately since an inequality that we
  use for the induction step applies only for $r>4$.  Let
  $H_1,H_2,\ldots,H_t$ be the large hyperplanes of $M$.  

  First consider $r=3$.  If $t\geq 3$ and $H_i\cap H_j=\emptyset$ for
  some $\{i,j\}\subsetneq [3]$, then we would get the contradiction
  $|H_1\cup H_2\cup H_3|\geq 2(2+g)+g=4+3g>|E(M)|$.  The same
  contradiction (with different counting) would occur if
  $H_1\cap H_2\cap H_3\ne\emptyset$.  Thus, $H_i\cap H_j$ is a
  singleton for all $\{i,j\}\subsetneq [3]$, while
  $H_1\cap H_2\cap H_3=\emptyset$.  Inclusion/exclusion now gives
  $|H_1\cup H_2\cup H_3|=3+3g$, so $H_1\cup H_2\cup H_3=E(M)$.  Thus,
  any large line contains at least two points in one of $H_1,H_2,H_3$
  and so is that large line, so $t=3$.  Also, the union of the three
  singletons $H_i\cap H_j$, for $\{i,j\}\subsetneq [3]$, is a basis,
  and the line spanned by any two of its elements is large.  Thus, (a)
  and (b) hold when $r=3$.

  Now consider $r=4$.  Focus on two large planes $H_1$ and $H_2$, and
  let $P$ be any other plane.  Since $|H_1|+|H_2|=6+6g=|E(M)|+2$, we
  have $|H_1\cap H_2|\geq 2$, so $|H_1\cap H_2|$ is $2$, $3$, or $2+g$
  by the assumption on the sizes of lines.  First assume that
  $|H_1\cap H_2|= 2$, so $H_1\cup H_2=E(M)$. Thus,
  $P=(P\cap H_1)\cup (P\cap H_2)$, so $|P|\leq 2(2+g)<3+3g$.  Thus,
  $H_1$ and $H_2$ are the only large planes.

  Now assume that $|H_1\cap H_2|= 3$.  Since
  $P\subseteq (P\cap H_1)\cup (P\cap H_2)\cup (E(M)-(H_1\cup H_2))$
  and $|E(M)-(H_1\cup H_2)|=1$, we get $|P|\leq 2(2+g)+1$.  If $g>2$,
  then $5+2g<3+3g$, and so $H_1$ and $H_2$ are the only large planes.
  Now assume that $g=2$, so large planes have nine elements and lines
  have at most four elements.  If $|P|=9$, then $P\cap H_1$ and
  $P\cap H_2$ are disjoint large lines, and so are disjoint from
  $H_1\cap H_2$.  Since $|H_1\cap H_2|= 3$, at most one large line of
  $H_1$ is disjoint from $H_1\cap H_2$, and likewise for $H_2$, so
  there is at most one large plane besides $H_1$ and $H_2$. Thus,
  there are at most three large planes.

  We may now assume that any two large planes of $M$ intersect in a
  large line.  No large line is in three large planes, for otherwise
  the union of the three large planes would contain $2+g+3(1+2g)=5+7g$
  elements, but $|E(M)|=4+6g$.  Since any large plane contains at most
  three large lines, each of which is in at most one more large plane,
  there are at most four large planes, so (a) holds when $r=4$.

  Assume that $M$ has four large planes, $H_i$ for $i\in[4]$. By the
  argument above, $H_i\cap H_j$ is a large line for all
  $\{i,j\}\subsetneq [4]$ and no two such intersections are equal, so
  each $H_i$ contains three large lines.  Using what we showed in the
  case $r=3$, it follows that $H_i\cap H_j\cap H_k$ is a singleton for
  all $\{i,j,k\}\subsetneq [4]$, all four such intersections are
  different, and their union is a basis, so (b) holds when $r=4$.

  Now assume that $r>4$ and that the bound on the number of large
  hyperplanes holds for lower ranks.  We claim that, if
  $ \{i,j\}\subseteq [t]$, then $H_i\cap H_j$ is a large rank-$(r-2)$
  flat.  Since $|H_i\cup H_j|\leq r+g\binom{r}{2}$, we have
  \begin{align*}
    |H_i\cap H_j|
    &\,= |H_i|+|H_j|-|H_i\cup H_j|\\
    &\,\geq 2\Bigl( r-1+g\binom{r-1}{2}\Bigr)- \Bigr(
      r+g\binom{r}{2}\Bigl)\\
    &\,= r-2+g\Bigl( 2\binom{r-1}{2}-\binom{r}{2}\Bigl)
    \\
    &\,= r-2+g\Bigl( \binom{r-2}{2}-1\Bigl)\\
    &\,> r-2 + g \binom{r-3}{2},
  \end{align*}
  (the final inequality requires $r>4$).  By the hypotheses,
  $|F|\leq r-3 + g \binom{r-3}{2}$ for a flat $F$ if $r(F) <r-2$, and
  $|F|\leq r-2 + g \binom{r-3}{2}$ if $r(F)=r-2$ but $F$ is not large,
  so $H_i\cap H_j$ is a large rank-$(r-2)$ flat, as claimed.

  Each of $H_i-(H_i\cap H_j)$ and $H_j-(H_i\cap H_j)$
  contains $$r-1+g\binom{r-1}{2}-(r-2) -g\binom{r-2}{2} =1+g(r-2)$$
  out of the $$r+g\binom{r}{2}-(r-2) -g\binom{r-2}{2} =2+g(2r-3)$$
  elements in $E(M)-(H_i\cap H_j)$, so $H_i$ and $H_j$ are the only
  large hyperplanes that contain $H_i\cap H_j$.  Therefore all
  $\binom{t}{2}$ intersections $H_i\cap H_j$, with
  $\{i,j\}\subseteq[t]$, are different.  In particular, $H_i\cap H_t$,
  for $i\in[t-1]$, are $t-1$ different large hyperplanes of $M|H_t$,
  so the induction assumption applied to $M|H_t$ gives $t-1 \leq r-1$,
  so $t\leq r$, so (a) holds.

  For the induction step in the proof that a basis
  $B=\{b_1,b_2,\ldots,b_r\}$ satisfying statement (b) exists, it will
  be useful to know some properties of such bases.  Counting shows
  that $\bigcup_{\{i,j\}\subsetneq[r]}\cl(\{b_i,b_j\})=E(M)$.  We
  claim that if $L$ is a line and $|L|\geq 3$, then we have
  $L\subseteq \cl(\{b_i,b_j,b_k\})$ for some $\{i,j,k\}\subseteq [r]$.
  Fix $\{e_1,e_2,e_3\}\subseteq L$.  There are $2$-subsets
  $B_1,B_2,B_3$ of $B$ for which $e_i\in\cl(B_i)$ for each $i\in [3]$.
  Now $e_3 \in \cl(B_1\cup B_2)\cap\cl(B_3)$, so
  $r(B_1\cup B_2\cup B_3)\leq 5$, so $B_1$, $B_2$, and $B_3$ are not
  pairwise disjoint, so we may assume that $B_1\cap B_2\ne\emptyset$.
  Thus, $L$ is contained in $\cl(B_1\cup B_2)$, which is either a line
  $\cl(\{b_i,b_j\})$ or a plane $\cl(\{b_i,b_j,b_k\})$, thus proving
  the claim.  It follows that the lines $\cl(\{b_i,b_j\})$ are the
  only large lines of $M$ since, for any large line $L$ of $M$, we
  have $L\subseteq \cl(\{b_i,b_j,b_k\})$ for some
  $\{i,j,k\}\subseteq [r]$, and the large plane $\cl(\{b_i,b_j,b_k\})$
  has three large lines, $\cl(\{b_i,b_j\})$, $\cl(\{b_i,b_k\})$, and
  $\cl(\{b_j,b_k\})$, as shown in the second paragraph of the proof.
  So, when a basis $B=\{b_1,b_2,\ldots,b_r\}$ as in (b) exists, it is
  unique.

  We can now prove that for any line $L$ with $b_i\in L$, either
  $|L|=2$ or $L=\cl(\{b_i,b_j\})$ for some $j\in[r]-\{i\}$.  This is
  because if $|L|\ne 2$, then $L\subsetneq \cl(\{b_i,b_j,b_k\})$ for
  some $\{j,k\}\subseteq [r]-\{i\}$, and
  $|L\cap \cl(\{b_j,b_k\})|\leq 1$, so $L$ has a point in one of
  $\cl(\{b_i,b_j\})$ or $\cl(\{b_i,b_k\})$ in addition to $b_i$, and
  so $L$ is that line.
  
  We now prove (b) by induction on $r$.  This was shown above for
  $r\in \{3, 4\}$.  For the induction step, assume that $r>4$ and that
  (b) holds for ranks less than $r$.  Thus, $M|H_1$ has a basis
  $\{b_1,b_2,\ldots,b_{r-1}\}$ that satisfies (b).  Now $M|H_2$ and
  $M|(H_1\cap H_2)$ also have such bases.  By the uniqueness of the
  basis of $M|(H_1\cap H_2)$ satisfying (b) and noting that the large
  lines of $M|(H_1\cap H_2)$ are large lines of $M|H_1$ and $M|H_2$ as
  well, it follows that, up to relabeling, the basis of $M|H_2$ that
  satisfies (b) is $\{b_2,b_3,\ldots,b_{r-1},b_r\}$ for some element
  $b_r$.  We claim that the basis of $M$ that satisfies (b) is
  $B=\{b_1,b_2,\ldots,b_r\}$.  To prove this, it suffices to show that
  the line $\cl(\{b_1,b_r\})$ is large since any other line spanned by
  two elements of $B$ is in at least one of $M|H_1$ and $M|H_2$ and so
  is large.  Consider the hyperplane $H=\cl(B-\{b_2\})$.  We know that
  the lines $\cl(\{b_i,b_j\})$, for $\{i,j\}\subsetneq [r]-\{2\}$ and
  $\{i,j\}\ne\{1,r\}$, are large since each is a large line of at
  least one of $M|H_1$ and $M|H_2$; thus,
  $|H|\geq r-1+g\bigl(\binom{r-1}{2}-1\bigl)$, so $H$ is a large
  hyperplane.  By the induction hypothesis, $M|H$ has a basis that
  satisfies (b), and, as we showed above, that basis accounts for all
  large lines of $M|H$.  Thus, that basis must be $B-\{b_2\}$, and so
  the line $\cl(\{b_1,b_r\})$ is large, as needed.  Thus, (b) holds.

  We now turn to (c).  As shown above, for the basis
  $B=\{b_1,b_2,\ldots,b_r\}$ in (b) and any $3$-point line $L$ of $M$,
  there is a subset $\{i,j,k\}$ of $[r]$ for which
  $L\subseteq \cl(\{b_i,b_j,b_k\})$; also, $L\cap B=\emptyset$.  The
  bound $\binom{r}{3}g^2$ follows since $\binom{r}{3}$ is the number
  of large planes $\cl(\{b_i,b_j,b_k\})$, and if
  $L\subseteq \cl(\{b_i,b_j,b_k\})$, then $L$ is determined by its
  intersections with any two of the large lines in
  $ \cl(\{b_i,b_j,b_k\})$, each of which has $g$ elements that are not
  in $B$.

  Now assume that $M$ has $\binom{r}{3}g^2$ three-point lines.
  Properties (1) and (2) of Theorem \ref{thm:axioms} follow from (b).
  Since $M$ has $\binom{r}{3}g^2$ three-point lines, each choice of a
  subset $\{i,j,k\}\subseteq[r]$ and elements
  $e\in \cl(\{b_i,b_j\})- \{b_i,b_j\}$ and
  $f\in\cl(\{b_i,b_k\})- \{b_i,b_k\}$ must give a $3$-point line
  $\cl(\{e,f\})$.  Thus, property (3) of Theorem \ref{thm:axioms} also
  holds, so, as claimed, $M$ is a Dowling geometry.
\end{proof}

The theorem and proof above easily yield the following alternative
characterization of Dowling geometries by optimizing valuative
invariants.  

\begin{thm}\label{thm:Dowling2}
  Fix integers $r\geq 3$ and $g\geq 2$.  Let $\mathcal{C}$ be the set
  of all simple rank-$r$ matroids on $r+g\binom{r}{2}$ elements in
  which lines have $2$, $3$, or $g+2$ points and, for each flat $F$ of
  rank $h$ with $2<h<r$, either $|F|=h+g\binom{h}{2}$ or
  $|F|\leq h+g\binom{h-1}{2}$.  Call rank-$i$ flats of size
  $i+g\binom{i}{2}$ \emph{large}, and call a flag \emph{large} if, for
  each $i\in[r]\cup\{0\}$, it has a large rank-$i$ flat.  Dowling
  geometries (including those over quasigroups if $r=3$) are the
  matroids in $\mathcal{C}$ that maximize the number of large flags
  and, after that, maximize the number of $3$-point lines.
\end{thm}

\begin{proof}
  Dowling geometries in $\mathcal{C}$ have $r!\,(g+2)/2$ large flags
  since a flag includes one of the $r$ large hyperplanes, and one of
  the $r-1$ large rank-$(r-2)$ flats that it contains, and so on
  through one of the three large lines in a large plane, and then one
  of the $g+2$ points in that large line.  Furthermore, since, by
  Theorem \ref{thm:Dowling}, matroids in $\mathcal{C}$ have at most
  $r$ large hyperplanes, each of which contains at most $r-1$ large
  rank-$(r-2)$ flats, etc., if a matroid in $\mathcal{C}$ has
  $r!\,(g+2)/2$ large flags, then it must have $r$ large hyperplanes.
  From that, part (c) of Theorem \ref{thm:Dowling} gives the result.
\end{proof}

\section{Spikes}

Spikes have played many important roles in matroid theory in recent
decades (see \cite{oxley}).  The following generalization of spikes
was introduced in \cite{genspike}: for $r, s, t\in\mathbb{N}$ with
$r\geq 3$, $s\geq r-1$, and $t\geq 2$, an \emph{$(r,s,t+1)$-spike with
  tip $a$} is a rank-$r$ simple matroid whose ground set is the union
of $s$ lines $L_1,L_2,\ldots,L_s$, each having $t+1$ elements and
containing $a$, for which, for any $k\in[r-1]$, the union of any $k$
of these lines has rank $k+1$.  An \emph{$r$-spike} is an
$(r,r,3)$-spike.  For example, $\PG(2,2)$ is a $3$-spike and
(atypically) any element can be the tip; for a given point $a$ of
$\PG(2,2)$, relaxing any one or more of the four circuit-hyperplanes
that do not contain $a$ gives more examples of $3$-spikes.  Also,
there is only one $(r,r-1,t+1)$-spike up to isomorphism, namely, the
parallel connection at $a$ of $r-1$ lines, each having $t+1$ elements.

In \cite{genspike}, $(r,s,t+1)$-spikes with $r\geq 5$ were shown to be
distinguishable from all other matroids (but not necessarily from each
other) by their Tutte polynomials, and it was shown which
$(r,s,t+1)$-spikes have the same Tutte polynomial; those results imply
that certain $(r,s,t+1)$-spikes are Tutte unique.  We identify some
$(r,s,t+1)$-spikes that are extremal.  (We treat only spikes with
tips.  Peter Nelson (private communications) has treated $r$-spikes
without tips.)

In any $(r,s,t+1)$-spike $M$ with $(t+1)$-point lines
$L_1,L_2,\ldots,L_s$ that contain $a$, any $3$-subset of any $L_i$ is
a circuit, as is any $4$-set $C$ for which, for some
$\{i,j\}\subseteq [s]$, we have
$|C\cap (L_i-\{a\})|=2=|C\cap (L_j-\{a\})|$.  All other circuits of
$M$ have size $r$ or $r+1$.  An $r$-set $X$ for which
$|X\cap(L_i-\{a\})|=1$ for $r$ integers $i\in[s]$ is either a basis or
a circuit; thus, there can be two types of $r$-circuits if
$r\in\{3,4\}$.  The hyperplane that is spanned by an $r$-circuit $X$
of the type just identified intersects each set $L_i-\{a\}$ in at most
one element, and if $s=r$, then $X$ is a circuit-hyperplane.

For any $s$, $t$, and $r>4$, there is a unique $(r,s,t+1)$-spike, up
to isomorphism, that has no $r$-circuits; this is the \emph{free
  $(r,s,t+1)$-spike}.  Similarly, the \emph{free $(3,s,t+1)$-spike} is
the $(3,s,t+1)$-spike in which the only $3$-circuits are the
$3$-subsets of the lines $L_i$, and the \emph{free $(4,s,t+1)$-spike}
is the $(4,s,t+1)$-spike in which the only $4$-circuits are those that
span a plane $L_i\cup L_j$.  Equivalently, the free $(r,s,t+1)$-spike
is the truncation to rank $r$ of the parallel connection of $s$ lines
at $a$, each having $t+1$ elements.  The free $(r,s,t+1)$-spike with
$r\geq 4$ is Tutte unique (see \cite{genspike}).

Among $3$-spikes, $\PG(2,2)$ maximizes the number of
circuit-hyperplanes.  For $r>3$, let $C$ be an $r$-circuit in
$\PG(r-1,2)$, let $a$ be an element of $\PG(r-1,2)$ that is not in the
hyperplane $\cl(C)$, and let $X$ be the union of the lines
$L_e=\cl(\{a,e\})$ for $e\in C$.  The restriction $M=\PG(r-1,2)|X$ is
an $r$-spike, and $C$ is a circuit-hyperplane of $M$.  It is easy to
show that the symmetric difference of $C$ with $k$ of the sets
$L_e-\{a\}$, for $e\in C$, is a circuit-hyperplane if and only if $k$
is even, otherwise the symmetric difference is a basis of $M$.  Thus,
$M$ has $2^{r-1}$ circuit-hyperplanes.  Also, for $r>3$, the maximum
number of circuit-hyperplanes in an $r$-spike is $2^{r-1}$, and using
the ideas about symmetric differences above, it is not hard to
construct an isomorphism between $M$ and any $r$-spike with $2^{r-1}$
circuit-hyperplanes.  This justifies calling the $r$-spike with
$2^{r-1}$ circuit-hyperplanes, for $r>3$, the \emph{binary $r$-spike}.
Of course, $\PG(2,2)$ is the binary $3$-spike. The binary $r$-spike
with $r\geq 5$ is Tutte unique (see \cite{genspike}).

The next result shows how to identify $(r,s,t+1)$-spikes, where
$r\geq 5$, by optimizing valuative invariants, and we show that the
free $(r,s,t+1)$-spike, the $(r,r-1,t+1)$-spike, and the binary
$r$-spike are extremal.  

\begin{thm}\label{thm:spikes}
  Fix integers $r\geq 5$, $t\geq 2$, and $s\geq r-1$. Let $M$ be a
  rank-$r$ matroid on $st+1$ elements.  Assume that for each
  $i\in[r-1]$, each rank-$i$ flat has at most $(i-1)t+1$ elements. For
  $i\in[r-1]$, call rank-$i$ flats with $(i-1)t+1$ elements
  \emph{large}.  Assume that $M$ minimizes the valuative invariant
  $N\mapsto F_{2,3}(N;t+1,w)$ for all $w\ne 2t+1$, so each plane that
  contains a large line is large.  Then $M$ has at most $s$ large
  lines, and if $M$ has $s$ large lines, then $M$ is an
  $(r,s,t+1)$-spike.  Also, the free $(r,s,t+1)$-spike, the binary
  $r$-spike, and the $(r,r-1,t+1)$-spike are extremal.
\end{thm}

\begin{proof}
  By the sizes of large lines and planes, coplanar large lines
  intersect in a point.  Let $L$ be any large line of $M$.  Any plane
  $P$ with $L\subsetneq P$ is large, so $|P-L|=t$.  Now
  $|E(M)-L|=(s-1)t$ and
  $\pi=\{P-L\,:\,P\text{ is a plane, }L\subsetneq P\}$ is a partition
  of $E(M)-L$, so $\pi$ has $s-1$ blocks; say
  $\pi=\{B_i\,:\,i\in[s-1]\}$.  The union of two planes $L\cup B_i$
  and $L\cup B_j$ has rank $4$ and size $3t+1$, so $L\cup B_i\cup B_j$
  is a flat. Let $L'$ be a line of $M$ for which $r(L\cup L')=4$.
  Thus, $\cl(L\cup L')=L\cup B_i\cup B_j$ for some
  $\{i,j\}\subseteq[s-1]$, and $|L'\cap B_i|=1= |L'\cap B_j|$, so
  $|L'|=2$.  Thus, any large line of $M$ other than $L$ is
  $B_i\cup \{a\}$ for some $a\in L$ and $i\in[s-1]$, so there are at
  most $s$ large lines.

  Since $L$ is any large line, the argument above shows that any two
  large lines of $M$ are coplanar.  All large lines $L'\ne L$ contain
  the same point of $L$ since two such large lines $L'$ and $L''$
  intersect in some point, which is therefore in the planes $L\cup L'$
  and $L\cup L''$ and so in their intersection, $L$.  Thus, if $M$ has
  $s$ large lines, then $M$ is an $(r,s,t+1)$-spike.

  The free $(r,s,t+1)$-spike is extremal since, among
  $(r,s,t+1)$-spikes, it is the only one that minimizes the number of
  $r$-circuits.  The binary $r$-spike is extremal since, among
  $r$-spikes, it is the only one that maximizes the number of
  circuit-hyperplanes.  Finally, the $(r,r-1,t+1)$-spike is extremal
  since it is the only spike with its parameters.
\end{proof}

We next give another characterization of $(r,r-1,t+1)$-spikes so that
in the next section we can show that direct sums of such spikes, where
$t$ is fixed, are extremal.  The notion of an $(r,r-1,t+1)$-spike
makes sense even when $r=2$, in which case a $(2,1,t+1)$-spike is a
$(t+1)$-point line, $U_{2,t+1}$.  What we gain by this mild extension
is a slightly stronger result about direct sums in the next section.
While Theorem \ref{thm:spikes} is limited to ranks $5$ and higher,
Theorem \ref{thm:spike2} applies for ranks $2$ and higher.

\begin{thm}\label{thm:spike2}
  Fix integers $r,t\in\mathbb{N}-\{1\}$.  Let $M$ be a rank-$r$
  matroid on $(r-1)t+1$ elements.  Assume that for each $i\in[r]$,
  each rank-$i$ flat has at most $(i-1)t+1$ elements. For $i\in[r]$,
  call rank-$i$ flats with $(i-1)t+1$ elements \emph{large}, and call
  a flag of flats of $M$ \emph{large} if it is counted by
  $F_{0,r}(N;0,1,t+1,2t+1,\ldots,(r-1)t+1)$.  Then $M$ has at most
  $(r-1)!(t+1)$ large flags, and $M$ has $(r-1)!(t+1)$ large flags if
  and only if $M$ is an $(r,r-1,t+1)$-spike.  Also, if $r>2$, then $M$
  has $(r-1)!(t+1)$ large flags if and only if $M$ has $r-1$ large
  hyperplanes, and that is the maximum number of large hyperplanes
  under the assumptions on $M$.
\end{thm}

\begin{proof}
  The statement is easy to see if $r=2$, so we focus on $r\geq 3$.  We
  first treat the bound on the number of large hyperplanes.  A
  rank-$3$ matroid on $2t+1$ points has at most two large lines, so
  the bound holds for $r=3$.  Now assume that $r>3$ and that the bound
  holds for smaller ranks.  Let $H_1,H_2,\ldots,H_k$ be the large
  hyperplanes of $M$.  For $\{i,j\}\subseteq[k]$, we have
  $H_j- H_i\subseteq E(M)-H_i$ and
  $|E(M)-H_i|=t\leq |H_j- (H_i\cap H_j)|$, so $E(M)-H_i=H_j- H_i$;
  likewise, $E(M)-H_j=H_i- H_j$.  Thus, $H_i\cup H_j= E(M)$ and
  $H_i\cap H_j$ is a large rank-$(r-2)$ flat.  Also, all intersections
  $H_i\cap H_j$, for $\{i,j\}\subseteq [k]$, are different since $H_i$
  and $H_j$ are the only hyperplanes that contain $H_i\cap H_j$.  The
  flats $H_i\cap H_k$, for $i\in[k-1]$, are large hyperplanes of
  $M|H_k$, so $k-1\leq r-2$ by the induction hypothesis, so
  $k\leq r-1$, as needed.

  A routine induction now shows that the intersection of any $i$ large
  hyperplanes is a large rank-$(r-i)$ flat.  Thus, all such
  intersections are different.  Also, if $k=r-1$, then
  $\bigcap_{i\in[r-1]}H_i$ is a rank-$1$ flat, say $\{a\}$, and $M$
  has $\binom{r-1}{r-2}=r-1$ large lines, all of which contain $a$, so
  $M$ is an $(r,r-1,t+1)$-spike, and $M$ has $(r-1)!(t+1)$ large
  flags.
  
  Arguments like those given in earlier sections show that, when
  $r>2$, having $r-1$ large hyperplanes is equivalent to having
  $(r-1)!(t+1)$ large flags.
\end{proof}
  
\section{Direct sums of some extremal matroids}\label{sec:dirsum}

Ferroni and Fink \cite{LuisAlex} conjectured that direct sums of
extremal matroids are extremal.  While it is far from settling the
conjecture, in Theorem \ref{thm:dirsumfamilies} below we show that
direct sums of extremal matroids are extremal if they share certain
features that can be described by optimizing valuative invariants.
The relevant common features occur in, for instance, (a) cycle
matroids of complete graphs, (b) projective geometries of rank at
least four and the same order, (c) affine geometries of rank at least
four and the same order, (d) Dowling geometries based on a group that
is the only group of its order, and (e) $(r,r-1,t+1)$-spikes where $t$
is fixed.
 
We start with another special case of the conjecture; we treat direct
sums of uniform matroids. Ferroni and Fink \cite[Example
4.4]{LuisAlex} observed that uniform matroids are extremal since they
maximize the number of bases; alternatively, Theorem
\ref{thm:PGviaFlags} applies since uniform matroids are perfect
matroid designs.  The proof of the next result uses the lexicographic
order on $(n,r)$-sequences: for two $(n,r)$-sequences
$\underline{x}=x_1\ldots x_n$ and $\underline{y}=y_1\ldots y_n$, we
have $\underline{x}<\underline{y}$ if, for the least $i$ with
$x_i\ne y_i$, we have $x_i<y_i$, so $x_i=0$ and $y_i=1$.
 
\begin{thm}\label{thm:dirsumunif}
  Direct sums of uniform matroids are extremal.
\end{thm}

\begin{proof}
  Let $M$ be
  $U_{r_1,n_1}\oplus U_{r_2,n_2}\oplus \cdots \oplus U_{r_k,n_k}$
  where $r_1\leq r_2\leq\cdots\leq r_k$ and, if $r_i=r_{i+1}$, then
  $n_i\geq n_{i+1}$.  Set $r=r(M)$ and $n=|E(M)|$.  If $M$ has $s$
  loops and $t$ coloops, then we can focus on matroids that have $s$
  loops and $t$ coloops by maximizing the valuative invariant
  $N\mapsto f_0(N;s,0)+ f_r(N;n,t)$.  Thus, without loss of
  generality, we can assume that $s=t=0$, so $0<r_i<n_i$ for all
  $i\in[k]$.  The sum of the coefficients of the terms $x^iy^j$, for
  $i+j<k$, in the Tutte polynomial of a matroid is a valuative
  invariant, and by minimizing it, we can focus on matroids that have
  at least $k$ connected components.  The number of cyclic flats in a
  matroid is a valuative invariant, and the matroids that minimize it,
  within the class of matroids that have at least $k$ connected
  components and neither loops nor coloops, have exactly $k$ connected
  components and their $2^k$ cyclic flats are the unions of none,
  some, or all of the connected components; therefore the restriction
  to such a connected component $X$ has just two cyclic flats,
  $\emptyset$ and $X$, so such restrictions are uniform.  Let
  $\underline{x}=x_1\ldots x_n$ be the $(n,r)$-sequence
  $$\underbrace{11\ldots 1}_{r_1\,\,1\text{'s}}\underbrace{0 0\dots
    0}_{n_1-r_1\,\,0\text{'s}} \underbrace{11\ldots
    1}_{r_2\,\,1\text{'s}}\underbrace{0 0\dots 0}_{n_2-r_2
    \,\,0\text{'s}} \ldots \underbrace{11\ldots
    1}_{r_k\,\,1\text{'s}}\underbrace{0 0\dots 0}_{n_k-r_k
    \,\,0\text{'s}}.
  $$ 
  Note that $\underline{x}$ encodes the rank and size of each
  connected component of $M$, and it is the lexicographically least
  rank sequence $\underline{r}$ for which $g_{\underline{r}}(M)$, the
  coefficient of $[\underline{r}]$ in $\mathcal{G}(M)$, is nonzero.
  For each $(n,r)$-sequence $\underline{y}$ with
  $\underline{y}<\underline{x}$ in lexicographic order, minimize the
  valuative invariant $N\mapsto g_{\underline{y}}(N)$.  Then maximize
  the valuative invariant $N\mapsto g_{\underline{x}}(N)$.  Let $N$ be
  a matroid that emerges from this sequence of optimization problems,
  so $N$ is an $n$-element rank-$r$ direct sum of $k$ uniform
  matroids for which $\underline{x}$ is the lexicographically least
  rank sequence $\underline{r}$ for which $g_{\underline{r}}(N)\ne 0$.
  Thus, the components of $N$ have the same ranks and size as those of
  $M$, so $N$ is isomorphic to $M$.  Thus, $M$ is extremal.
\end{proof}

We next give a result on direct sums of extremal matroids that applies
to direct sums of many of the matroids that we treated earlier in this
paper.  Looking ahead to Corollary \ref{cor:extKnPGAGD} may make the
notation below easier to digest.  Characterizations of the type used
in this result are illustrated by Theorem \ref{thm:Kn2} for
$M(K_{r+1})$.  In many applications, $\mathcal{M}$ contains, up to
isomorphism, one matroid of each rank (e.g,. $\mathcal{M}$ can be the
set of cycle matroids of complete graphs of rank two and greater).
For the general case, an example to keep in mind is Dowling
geometries, where many matroids maximize the number of large flags,
but only Dowling geometries also maximize the number of $3$-point
lines.

\begin{thm}\label{thm:dirsumfamilies}
  Let $\alpha_0=0,\alpha_1=1,\alpha_2,\alpha_3,\ldots$ be integers for
  which 
  \begin{itemize}
  \item[(1)] $\alpha_2\geq 3$ and
    $\alpha_i+\alpha_{j}\leq \alpha_{i+j-1}+1$ for all
    $i,j\in\mathbb{N}$.
  \end{itemize}
  Fix a set $A\subseteq [\alpha_2]-\{1\}$ with $2,\alpha_2\in A$. Fix
  an integer $d\geq 2$.  Let $\mathcal{M}$ be a set of simple matroids
  of rank at least $d$ that satisfies the following conditions:
  \begin{itemize}
  \item[(2)]
    for all $M\in\mathcal{M}$ and $i\leq r(M)$, some rank-$i$
    flat of $M$ has size $\alpha_i$,
  \item[(3)] for each $M\in\mathcal{M}$ and flat $F$ of $M$, if
    $2<r(F)<r(M)$, then either $|F|= \alpha_{r(F)}$ or
    $|F|\leq \alpha_{r(F)-1}+1$, while if $r(F)=2$, then $|F|\in A$,
    and
  \item[(4)] for $r\geq d$, a rank-$r$ simple matroid $M$ is in
    $\mathcal{M}$ if and only if, among rank-$r$ simple matroids that
    satisfy conditions \emph{(2)-(3)}, the valuative invariant that
    maps $N$ to
    $F_{0,r}(N;\alpha_0,\alpha_1,\alpha_2,\ldots,\alpha_r)$ is
    maximized by $M$.
  \end{itemize}
  If, for some sequence of valuative invariants $f_1,f_2,\ldots,f_k$
  that are additive on direct sums (i.e.,
  $f_i(M\oplus N)=f_i(M)+f_i(N)$), the matroids $M_1,M_2,\ldots,M_t$
  are, up to isomorphism, the only matroids in $\mathcal{M}$ of their
  ranks that maximize $f_1$, and then $f_2$, and so on, through $f_k$,
  then $M_1\oplus M_2\oplus \cdots\oplus M_t$ is extremal.
\end{thm}

\begin{proof}
  Call rank-$i$ flats of size $\alpha_i$ in a matroid
  $M\in\mathcal{M}$ \emph{large}, and call a flag of $M$ \emph{large}
  if it has flats of each rank from $0$ to $r(M)$, all of which are
  large.  Thus, if $M\in \mathcal{M}$, then
  $F_{0,r}(M;\alpha_0,\alpha_1,\alpha_2,\ldots,\alpha_r)$ is the
  number of large flags of $M$.  It follows from conditions (1)-(2)
  that each matroid in $\mathcal{M}$ is connected.
  
  Fix any matroids $M_1,M_2,\ldots,M_t$ in $\mathcal{M}$ on disjoint
  ground sets.  We can assume that $r(M_i)\geq r(M_{i+1})$ for
  $i\in[t-1]$.  Set $M=M_1\oplus M_2\oplus \cdots\oplus M_t$.  We
  start with some observations about $M$.

  By condition (1), for $i\geq 2$, no rank-$i$ flat in $M$ that has
  elements in two or more of $M_1,M_2,\ldots,M_t$ has size more than
  $\alpha_{i-1}+1$.  Thus, $M$ satisfies the bounds in condition (3)
  and any rank-$i$ flat of $M$ of size $\alpha_i$ is a flat of some
  $M_j$.  Therefore, we can focus on direct sums that satisfy the
  bounds in condition (3).

  For $s\in[t+1]$, let $p_s$ be the sum of the first $s-1$ terms among
  $r(M_1),r(M_2),\ldots,r(M_t)$, so $p_1=0$ and $p_{t+1}=r(M)$.  If
  $i\in [r(M)-1]$, then $i=p_s+u$ for unique $s\in[t]$ and
  $u\in[r(M_s)-1]\cup\{0\}$, and the largest rank-$i$ flats in $M$
  have size
  $$\beta_i=\alpha_{r(M_1)}+\alpha_{r(M_2)}+\cdots
  +\alpha_{r(M_{s-1})}+\alpha_u.$$ Set $\beta_{r(M)}=|E(M)|$.  We
  extend the terminology of large flats and large flags to $M$ in the
  obvious way, using the sequence
  $\beta_0,\beta_1,\ldots,\beta_{r(M)}$.  Note that
  $\beta_i\leq\alpha_i$ for all $i$, with equality if and only if
  $i\leq r(M_1)$.  To get the large flags of $M$, start with a large
  flag of a matroid among $M_1,M_2,\ldots,M_t$ of largest rank, then
  take the union of the ground set of that matroid with, in turn, the
  nonempty sets in a large flag of a matroid of largest rank among
  $M_1,M_2,\ldots,M_t$ that has not yet been selected, and repeat
  until we reach $E(M)$.  Condition (1) gives $\beta_{k+1}-\beta_k=1$
  for $t$ integers $k\in[r(M)]$, namely, $k=p_s$ with $s\in[t]$.  If,
  for some $s\in[t]$, a sequence
  $\gamma_0,\gamma_1,\ldots,\gamma_{r(M)}$ has
  $\gamma_{p_s}=\beta_{p_s}$ and $\gamma_{p_s+1}-\gamma_{p_s}>1$, then
  $\gamma_{p_s+1}>\beta_{p_s+1}$ and so
  $F_{0,r(M)}(M;\gamma_0,\gamma_1,\ldots,\gamma_{r(M)})=0$ since
  $\beta_{p_s+1}$ is the largest size of a rank-$(p_s+1)$ flat of $M$.

  Consider matroids of rank $r(M)$ and size $ |E(M)|=\beta_{r(M)}$.
  We claim that optimizing the valuative invariants below, in the
  order listed, isolates the direct sums
  $N_1\oplus N_2\oplus \cdots\oplus N_t$ where $N_i\in\mathcal{M}$
  and, if the components are indexed so that $r(N_i)\geq r(N_{i+1})$
  for all $i\in[t-1]$, then $r(M_i)=r(N_i)$ for all $i\in[t]$.
  \begin{itemize}
  \item[(a)] To focus on simple matroids, maximize
    $N\mapsto f_1(N;1)$.
  \item[(b)] As justified above, minimize $N\mapsto f_2(N;j)$ for each
    $j\not \in A$; also, for all $i, j$ with $3\leq i\leq r(M)$ and
    $\alpha_{i-1}+1< j<\alpha_i$, minimize $N\mapsto f_i(N;j)$.
  \item[(c)] To focus on matroids with at least $t$ connected
    components, minimize the sum of the coefficients of the terms
    $x^iy^j$ with $i+j<t$ in the Tutte polynomial $T(N;x,y)$.
  \item[(d)] For each sequence
    $\gamma_0,\gamma_1,\ldots,\gamma_{r(M)}$ where, for some
    $s\in[t]$, we have $\gamma_{p_s}=\beta_{p_s}$ and
    $\gamma_{p_s+1}-\gamma_{p_s}>1$, minimize
    $N\mapsto F_{0,r(M)}(N;\gamma_0,\gamma_1,\ldots,\gamma_{r(M)})$.
  \item[(e)] Maximize
    $N\mapsto F_{0,r(M)}(N;\beta_0,\beta_1,\ldots,\beta_{r(M)})$.
  \end{itemize}
  Since $\beta_{k+1}-\beta_k=1$ if and only if $k=p_s$ where
  $s\in[t]$, simple matroids that satisfy (c) and (e) have $t$
  connected components.  Let $N=N_1\oplus N_2\oplus\cdots\oplus N_t$
  satisfy the sequence of conditions (a)-(e), so each $N_i$ is
  connected.  Let
  $F_0\subsetneq F_1\subsetneq \cdots\subsetneq F_{r(M)}$ be a flag
  that $F_{0,r(M)}(N;\beta_0,\beta_1,\ldots,\beta_{r(M)})$ counts.
  Since $\beta_{r(M_1)}=\alpha_{r(M_1)}$, we have
  $F_{r(M_1)}\subseteq E(N_j)$ for some $j\in[t]$; we may assume that
  $j=1$.  If $ E(N_1)\ne F_{r(M_1)}$, then $|X-F_{r(M_1)}|>1$ for some
  flat $X$ of $N_1$ that covers $F_{r(M_1)}$ (otherwise more than $t$
  sets $F_{i+1}-F_i$ would be singletons, so more than $t$ terms
  $\beta_{i+1}-\beta_i$ would be $1$); thus,
  $F_{0,r(M)}(N;\gamma_0,\gamma_1,\ldots,\gamma_{r(M)})>0$ for some
  sequence $\gamma_0,\gamma_1, \ldots,\gamma_{r(M)}$ where
  $\gamma_i=\beta_i$ for $i\leq r(M_1)$ and
  $\gamma_{r(M_1)+1}-\gamma_{r(M_1)}=|X-F_{r(M_1)}|>1$, but the
  minimization in item (d) prevents that.  Thus, $E(N_1)= F_{r(M_1)}$.
  The same ideas show that, up to reindexing, $r(N_i)=r(M_i)$ and
  $E(N_i)=F_{p_{i+1}}- F_{p_i}$ for all $i\in [t]$.  We maximize
  $F_{0,r(M)}(N;\beta_0,\beta_1,\ldots,\beta_{r(M)})$ when each
  $ F_{0,r(N_i)}(N_i;\alpha_0,\alpha_1,\ldots,\alpha_{r(N_i)})$ is
  maximized, so each $N_i$ is in $\mathcal{M}$, as needed.

  We now use the hypothesis that there are valuative invariants
  $f_1,f_2,\ldots,f_k$ for which, among all matroids in $\mathcal{M}$,
  the matroids $M_1,M_2,\ldots,M_t$ are, up to isomorphism, the unique
  matroids of their ranks that maximize $f_1$, and then $f_2$, and so
  on, through $f_k$. By the assumption that each $f_i$ is additive on
  direct sums, $f_i(N)$ is maximized precisely when each $f_i(N_j)$ is
  maximized.  By assumption, that is when $N_j$ is isomorphic to
  $M_j$.  Thus, maximizing $f_1$, and then $f_2$, and so on through
  $f_k$ leaves one matroid, up to isomorphism, namely, $M$, as
  claimed.
\end{proof}

The assumption that $f_1,f_2,\ldots,f_k$ are additive can be replaced
by any condition that ensures that each $f_i$ is maximized on a
matroid precisely when it is maximized on the restrictions to its
connected components.

\begin{cor}\label{cor:extKnPGAGD}
  Direct sums of cycle matroids of complete graphs of rank at least
  two are extremal.  For a fixed prime power $q$, direct sums of
  projective geometries of order $q$ and rank at least four are
  extremal, and the same holds for affine geometries.  If there is
  only one projective or affine plane of order $q$, up to isomorphism,
  then planes and lines of order $q$ can be included in the previous
  statement.  For a group $G$ for which $G$ is the only group of order
  $|G|$, direct sums of Dowling geometries of rank at least four based
  on $G$ are extremal.  For a fixed integer $t\geq 2$, direct sums of
  $(r,r-1,t)$-spikes of rank at least two are extremal.
\end{cor}

\begin{proof}
  For all cases except Dowling geometries, the set $\mathcal{M}$ in
  Theorem \ref{thm:dirsumfamilies} has, up to isomorphism, a unique
  matroid of each rank.  The first statement follows directly from
  Theorems \ref{thm:Kn2} and \ref{thm:dirsumfamilies} with
  $\alpha_i=\binom{i+1}{2}$, $A=\{2,3\}$, and $d=2$.  The assertion
  about projective geometries follows Theorems \ref{thm:PGviaFlags}
  and \ref{thm:dirsumfamilies} with $\alpha_i=\frac{q^i-1}{q-1}$,
  $A=\{2,q+1\}$, and $d=4$; likewise for affine geometries, with
  $\alpha_i=q^{i-1}$, for $i\geq 1$, $A=\{2,q\}$, and $d=4$; when, as
  noted above, those results extend to planes and lines, use $d=2$.
  For Dowling geometries, the result follows from Theorems
  \ref{thm:Dowling2} and \ref{thm:dirsumfamilies} with
  $\alpha_i=i+|G|\,\binom{i}{2}$, $A=\{2,3,|G|+2\}$, and $d=4$ by
  using a single additional valuative invariant $f_1$, namely, the
  number of $3$-point lines.  For $(r,r-1,t)$-spikes, the result
  follows from Theorems \ref{thm:spike2} and \ref{thm:dirsumfamilies}
  with $\alpha_i=1+(i-1)t$ for $i\geq 1$, $A=\{2,t+1\}$, and $d=2$.
\end{proof}

This result specifies rank at least four for Dowling geometries since
only for the first few orders $g$ is there a unique quasigroup of
order $g$.  That stands in contrast to groups, where there is a unique
group of order $g$ if and only if $\gcd(g,\phi(g))=1$, where $\phi(g)$
is Euler's $\phi$-function; see \cite{uniquegp}.

\medskip

\begin{samepage}

\begin{center}
 \textsc{Acknowledgments}
\end{center}

\vspace{3pt}

The author is grateful to Luis Ferroni, Alex Fink, and Peter Nelson
for helpful and encouraging e-mail exchanges while this material was
being developed.  The author thanks the referee for valuable comments
and especially for mentioning Lemma \ref{lem:refobs}.

\end{samepage}

\end{document}